\documentclass[a4paper,12pt,oneside]{article}
\usepackage[latin1]{inputenc}
\usepackage[sumlimits,namelimits]{amsmath}
\usepackage{amssymb,amscd,latexsym}
\usepackage{theorem}
\usepackage{enumerate}
\usepackage[T1]{fontenc}
\usepackage{graphicx}
\usepackage{multicol}
\usepackage{pstcol,pst-char}

\newtheorem{dfn}{Definition}[section]
\newtheorem{prp}[dfn]{Proposition}
\newtheorem{thm}[dfn]{Theorem}

\newtheorem{lmm}[dfn]{Lemma}
\newtheorem{rmk}[dfn]{Remark}

\def\eps{\varepsilon}
\def\R{\mathbb{R}}
\def\N{\mathbb{N}}
\def\S{\mathcal{S}}
\def\D{\mathcal{D}}
\def\B{\mathcal{B}}
\def\M{\mathcal{M}(m,\kappa,r_0)}
\def\E{\mathcal{E}(n,k_1,k_2)}
\def\FnX{\mathcal{F}(X)}
\def\FnXX{\mathcal{F}(X \times X)}

\numberwithin{equation}{section}

\title{Discretization of vector bundles and rough Laplacian}
\author{Tatiana Mantuano
\thanks{Supported by Swiss National Science Foundation, grant No. 20-101469.}
}
\date{}

\begin{document}

\maketitle

{\abstract Let $\M$ be the set of all compact connected
$m$-dimensional manifolds $(M,g)$ such that $Ricci(M,g) \geq
-(m-1) \kappa g$ and $Inj(M,g) \geq r_0>0$. Let $\E$ be the set of
all Riemannian vector bundles $(E,\nabla)$  of real rank $n$ with
$|R^E| \leq k_1$ and $|d^* R^E| \leq k_2$. For any vector bundle
$E \in \E$ with harmonic curvature or with complex rank one, over
any $M \in \M$ and for any discretization $X$ of $M$ of mesh $0
<\eps \leq \frac{1}{20}r_0$, we construct a canonical twisted
Laplacian $\Delta_A$ and a potential $V$ depending only on the
local geometry of $E$ and $M$ such that we can compare uniformly
the spectrum of the rough Laplacian $\overline{\Delta}$ associated
to the connection of $E$ and the spectrum of $\Delta_A + V$. We
show that there exist constants $c$, $c' > 0$ depending only on
the parameters of $\M$ and $\E$ such that $c'\lambda_k(X,A,V) \leq
\lambda_k(E)\leq c\lambda_k(X,A,V)$, where $\lambda_k(\cdot)$
denotes the $k^{th}$ eigenvalue of the considered operators ($k
\leq n |X|$). For flat vector bundles, we show that the potential
is zero, $\Delta_A$ turns out to be a discrete magnetic Laplacian
and we relate $\lambda_1(E)$ to the holonomy of $E$.

\smallskip

\textbf{Mathematics Subject Classification (2000)}: 58J50, 53C20.

\smallskip

\textbf{Key words}: connection, rough Laplacian, discrete magnetic
La\-placian, Harper operator, eigenvalues, discretization,
holonomy.

}

%\tableofcontents

\section{Introduction}

In \cite{Ma}, we have shown  that for a family of compact
connected manifolds $\M$ with injectivity radius and Ricci
curvature bounded below (i.e. $(M,g)\in \M$ if $M$ is a compact
connected $m$-dimensional Riemannian manifold with $Ricci(M,g)
\geq -(m-1)\kappa g$ and $Inj(M,g) \geq r_0$), we can compare
uniformly the spectrum of the Laplacian acting on functions with
the spectrum of the combinatorial Laplacian acting on a graph with
fixed mesh constructed on the manifolds. Indeed, we show that
there exist positive constants $c$, $c'$ depending on the
parameters of the problem such that for any $M \in \M$ and any
discretization $X$ of $M$ (with mesh $\eps < \frac{1}{2}r_0$), the
following holds
\begin{eqnarray}
    c' \lambda_k(X)\leq \lambda_k(M) \leq c\lambda_k(X) \label{result1}
\end{eqnarray}
for $k < |X|$, where $\lambda_k(\cdot)$ stands for the $k^{th}$
eigenvalue of the considered Laplacian. This result generalizes in
a natural way different works like \cite{Br1}, \cite{Bur},
\cite{Bu} and \cite{Ka2} that were motivated either by the study
of the relation between the fundamental group of a manifold and
the spectrum of its finite coverings (\cite{Br1}, \cite{Bur}) or
by the relation between the spectrum of a manifold and its Cheeger
isoperimetric constant  (\cite{Bu}) or by the existence of
harmonic functions (\cite{Ka2}). More generally, the aim of the
discretization is to have an understanding of the spectrum (a
global invariant on the manifold) with a minimum of informations
about the local geometry of the manifold.

\medskip

Of course, the problem is interesting for differential operators
other than the Laplacian and we may address the following
question: does the same kind of comparison hold for other
geometric differential operators such that the Laplacian acting on
$p$-forms or the Dirac operator? Most of these operators may be
expressed in terms of a connection Laplacian added with a
curvature term. In this article, we investigate the case of such a
connection (or rough) Laplacian $\overline{\Delta}$ associated to
a connection $\nabla$ on a vector bundle. More precisely, the
purpose is to establish a uniform comparison of spectra between
rough Laplacians on vector bundles and twisted Laplacians on
graphs that generalize combinatorial or discrete magnetic
Laplacians. The Riemannian vector bundles we are interested in
have curvature and exterior coderivative of curvature bounded i.e.
we study Riemannian vector bundles $E$ with fiber of real rank $n$
such that $|R^E | \leq k_1$ and $|d^* R^E| \leq k_2$ (denote by
$\E$ the set of such vector bundles). The main result (Theorem
\ref{main}) states that there exist positive constants $c$, $c'$
(depending only on the given parameters) such that for any vector
bundle $E \in \E$ over any $M \in \M$ satisfying one of the
following assumptions
\begin{enumerate}
    \item[I)] the curvature of $E$ is harmonic i.e. $d^* R^E =0$,
    \item[II)] $E$ is of complex (or quaternionic) rank one
\end{enumerate}
and for any discretization $X$ of $E$, we can construct a
canonical twisted Laplacian $\Delta_A$ and a potential $V$
depending only on the local geometry of $E$ such that
\begin{eqnarray}
    c' \lambda_k(X,A,V)\leq \lambda_k(E) \leq c\lambda_k(X,A,V) \label{result2}
\end{eqnarray}
for any $k \leq n |X|$, where $\lambda_k(E)$ denotes the $k^{th}$
eigenvalue of the rough Laplacian $\overline{\Delta}$ and
$\lambda_k(X,A,V)$ the $k^{th}$ eigenvalue of $\Delta_A + V$.

\medskip

The case of flat vector bundles is especially enlightening.
Indeed, if $E$ is flat, we show that the potential $V$ is zero and
that $\Delta_A$ is a discrete magnetic Laplacian. This particular
case shows how the construction of $\Delta_A$ is strongly related
to the holonomy of $E$. This fact is emphasized by Theorem
\ref{thmHolo} which relates the holonomy (in the sense of
\cite{BBC}) to the first eigenvalue of $\Delta_A$ and therefore of
$\overline{\Delta}$. In order to understand the problem of
non-flat vector bundles, go back to the case of functions. Recall
that for functions we needed to establish correspondances between
functions on the manifold and functions on the graph. To that aim
and in particular to associate smooth functions to functions on
the graph, we had to extend locally such a function in a constant
way and then smooth it (with a partition of unity). The question
of extending locally is a central problem for the case of vector
bundles. It turns out that extending by parallel transport is
really efficient for flat vector bundles as it produces parallel
sections. But, as soon as the curvature comes in, parallel
transport is not convenient anymore and we need to construct a
finer way to extend locally a section. In fact, the obstruction to
extend in a parallel manner is double: the holonomy plays the role
of a global obstruction to extend as parallel as possible and
locally the curvature plays the same role. The twisted Laplacian
will precisely render the holonomy of the vector bundle, while the
potential will take into account the local non-flat geometry.

\medskip

The paper is organized as follows.  In Section
\ref{sectionSettings}, we introduce the notations, we define the
general notion of twisted Laplacian on a graph and recall the main
properties of the discretization of a manifold (that will coincide
with the notion of discretization of vector bundles). Section
\ref{comp} is devoted to the proof of the main result (Theorem
\ref{main}). The main difficulty is to construct a suitable
twisted Laplacian (see Section \ref{Exten}). From a geometric
point of view, the problem is the dependence on the local geometry
of the Laplacian and the potential to have enough informations to
estimate the spectrum of the vector bundle. Technically, we need
fine analysis on vector bundles like Sobolev inequalities for
sections to achieve the construction. The particular case of flat
vector bundles can be kept in mind as the ground example during
the reading. In this case, the proofs can be done easier (we can
avoid the technical tools described in Section \ref{Exten}).
Nevertheless, this case already contains the essential information
for $\Delta_A$ as it shows how the holonomy is related to
$\Delta_A$ (see Section \ref{holonomy}). For non-flat vector
bundles, $\Delta_A$ does not suffice anymore to control the rough
Laplacian, so that we have to introduce a potential $V$ which
takes care of the curvature locally. The generalization of the
flat case is then done for two different cases (see assumptions
$\mathrm{I})$ and $\mathrm{II})$), for rank one vector bundles and
for vector bundles with harmonic curvature. These two cases are
really of different nature. This appears all along Section
\ref{comp} and this begins with the construction of $\Delta_A +V$
(in Section \ref{Constr}) which differs according to the
assumptions $\mathrm{I})$ or $\mathrm{II})$. In Section
\ref{holonomy}, we establish the relationship between the holonomy
and the first eigenvalue of the rough Laplacian for flat vector
bundles. The part of Theorem \ref{thmHolo} that bounds from below
the first eigenvalue in terms of the holonomy can be generalized
easily to vector bundles with harmonic curvature. But this will
not be done here. This result is in fact due to Ballmann, Brüning
and Carron in a more general setting (see \cite{BBC}). Finally, we
collect some more technical proofs in the appendix to make easier
the reading, even if the results are not of minor importance for
the paper.

\section{Settings}\label{sectionSettings}

\subsection{Rough Laplacian}

In this section, we recall basic facts on the rough Laplacian (for
a general reference see \cite{Ber}, \cite{Mo} or \cite{Pe} for
instance). Let $(M,g)$ be a compact connected $m$-dimensional
Riemannian manifold without boundary and with volume form denoted
by $dV$. Moreover, let $(E, \nabla)$ be a Riemannian vector bundle
with $n$-dimensional fiber over $M$ i.e. $E$ is a vector bundle
over $M$ endowed with a smooth metric $\langle \cdot , \cdot
\rangle$ and a compatible connection $\nabla$. On the set
$\Gamma(E)$ of smooth sections of $E$, denote by $(\cdot,\cdot)$
the $L^2$-inner product  endowed by $\langle \cdot , \cdot
\rangle$ and $g$. Recall that the connection extends to
$p$-tensors on $M$ with values in $E$ and that we define
$\nabla^{\ast}$ to be the adjoint of $\nabla$ with respect to the
$L^2$-inner product. The rough Laplacian (or connection Laplacian)
acting on $\Gamma(E)$ is then defined by $ \overline{\Delta} =
\nabla^{\ast} \nabla $. The spectrum of $\overline{\Delta}$ is
discrete and non-negative and will be denoted
\begin{eqnarray}
    Spec(E) = \{\lambda_1(E) \leq \lambda_2(E)\leq \ldots \leq
    \lambda_k(E) \leq \ldots\}\nonumber.
\end{eqnarray}
The Rayleigh quotient of a non-zero section $s$ is defined by
 $   R(s) = \frac{\|\nabla s\|^2}{\|s\|^2} $,
where $\|\cdot\|$ denotes the $L^2$-norm associated to the
$L^2$-inner product defined above. Later we will need the
following variational characterizations of $Spec(E)$ known as
min-max and max-min theorems. For any $k \geq 1$,
\begin{eqnarray}
    \lambda_k (E)& =&
    \min_{\Omega^k}\max \{R(s) : s \in \Omega^k\setminus\{0\}\}
    \nonumber\\
    &=&
    \max_{\Omega^{k-1}}\min \{R(s) : s \in \Omega^{k-1}\setminus\{ 0 \} ,
    \ s \bot \Omega^{k-1}\text{ w.r.t } (\cdot , \cdot ) \}
    \nonumber
\end{eqnarray}
where $\Omega^k$ (resp. $\Omega^{k-1}$) ranges over all
$k$-dimensional (resp. $(k-1)$-dim.) subspaces of $\Gamma(E)$.

%%%%%%%%%%%%%%%%%%%%%%%%%%%%%%%%%%%%%%%%%%%%%%%%%%%%%%%%%%%%%%%%%%%%%%%%%%%%%%%%%

\subsection{Twisted Laplacian}

Let $\Gamma = (X,E(X))$ be a finite connected graph endowed with
the path metric. For $p \in X$ denote by $N(p)$ the set of
vertices at distance $1$ from $p$ and by $m(p)$ the number of such
vertices. In order to generalize the combinatorial Laplacian (see
\cite{Lu} for a definition) and the discrete magnetic Laplacian
(see \cite{MY} for a definition), let us consider the set of
functions on $X$ with values in $\mathbb{R}^n$ i.e. $\FnX = \{ f :
X \rightarrow \mathbb{R}^{n} \}$, provided with the inner product
$ (f,g) = \sum_{p \in X} f(p) \cdot g(p)$, where $\cdot$ denotes
the Euclidean inner product of $\mathbb{R}^n$.
\begin{dfn} \label{dfnLapA}
    For any $p \in X$ and $q \in N(p)$ assume that $A(p,q) : \R^n\to \R^n$
    is a given linear transformation. The \textbf{twisted Laplacian}
    associated to $A$ is the operator $\Delta_A : \FnX \to \FnX$ defined by
\begin{eqnarray}
    \Delta_A f(p) = \frac{1}{2}\sum_{q \in N(p)}
    \left(\mathbb{I} + A^t(p,q) A(p,q)\right)f(p) - \left(
    A(q,p) + A^t(p,q) \right)f(q). \nonumber
\end{eqnarray}
\end{dfn}
\begin{rmk}
    If for any $ p $, $q$, the operator $A(p,q)$ is the identity,
    then $\Delta_A$ is the combinatorial Laplacian.
\end{rmk}
\begin{rmk}
    If $A(p,q)$ belongs to $O(n)$ and $A^t(p,q) = A(q,p)$, then
    \linebreak
    $\Delta_A f(p) = m(p) f(p) - \sum_{q \in N(p)} A(q,p)f(q)$. In
    this case the twisted Laplacian
    is usually called \textbf{discrete magnetic Laplacian} or Laplacian associated to the
    Harper operator $A$.
\end{rmk}
Let us introduce the space of functions $\FnXX = \{ F : X \times X
\to \mathbb{R}^n \}$ and provide it with the inner product given
by $ (F,G) = \frac{1}{2}\sum\limits_{p \in X}\sum\limits_{q \in X}
F(p,q)\cdot G(p,q)$.
\begin{lmm}
    Let $A(p,q)$ be as in Definition \ref{dfnLapA} and $\Delta_A$ the
    twisted Laplacian associated to $A$. Let $D_A : \FnX \to \FnXX$ be defined
    by
\begin{eqnarray}
    D_A f(p,q) = \left\{ \begin{array}{ll}
                            f(q)-A(p,q)f(p) & if \ p \in N(q),\\
                            0               & otherwise .
                          \end{array}  \right. \nonumber
\end{eqnarray}
Then, for any $f$,$g \in \FnX$, we have $(\Delta_A f,g) = (D_A f,
D_A g)$.
\end{lmm}
\textbf{Proof}: let $f$, $g \in \FnX$. Then, we have
\begin{eqnarray}
    (\Delta_A f , g)    & = &
    \frac{1}{2} \sum_{p \in X}\sum_{q \in N(p)}
    \left(f(p)-A(q,p)f(q)\right)\cdot g(p) \nonumber\\
    & & -
    \frac{1}{2} \sum_{p \in X}\sum_{q \in N(p)}
    \left(f(q)-A(p,q)f(p)\right)\cdot A(p,q)g(p)\nonumber\\
                        & = &
    \frac{1}{2} \sum_{p \in X}\sum_{q \in X} D_A f(q,p)\cdot g(p)\nonumber +
    \frac{1}{2} \sum_{p \in X}\sum_{q \in X} D_A f(p,q)\cdot D_A g(p,q)\nonumber\\
    &&-
    \frac{1}{2} \sum_{p \in X}\sum_{q \in X} D_A f(p,q)\cdot g(q)\nonumber
                        \ = \ (D_A f , D_A g).\nonumber \text{ \hspace{2cm} } \Box
\end{eqnarray}
A direct consequence of this lemma is that $\Delta_A$ is symmetric
and non-negative, so it admits a non-negative spectrum. If $V :
\FnX \to \FnX$ is a non-negative potential, then the spectrum of
$\Delta_A + V$ is characterized by min-max theorem as follows
\begin{eqnarray}
    \forall 1 \leq k \leq n |X|, \ \ \ \ \
    \lambda_k(X,A,V) = \min\limits_{W^k} \max \{ R(f) : f \in W^k
    \setminus\{0\}\} \nonumber
\end{eqnarray}
where $W^k$ ranges over all $k$-dimensional vector subspaces of
$\FnX$ and $R(f)$ is the Rayleigh quotient of $f$ defined by $
R(f) = \frac{\|D_A f\|^2+ (Vf,f)}{\|f\|^2}$.

%%%%%%%%%%%%%%%%%%%%%%%%%%%%%%%%%%%%%%%%%%%%%%%%%%%%%%%%%%%%%%%%%%%%%%%%%%%%%%%%%%

\subsection{Discretization of vector bundles}

In this section, we define the notion of discretization of a
vector bundle.

\begin{dfn}
Let $(E, \nabla)$ be a Riemannian vector bundle over $(M,g)$ a
compact connected Riemannian manifold with $\partial M =
\emptyset$. An \textbf{$\eps$-discre\-ti\-zation} of $E$ is a
discre\-ti\-zation of $M$ of mesh $\eps >0$.
\end{dfn}
The discretization of a manifold (of mesh $\eps$) is defined as in
\cite{Ch} (Section V.3.2). Let us recall the definition and the
properties of such a discretization. Let $(M , g)$ be a compact
connected $m$-dimensional Riemannian ma\-nifold. A discretization
of $M$, of mesh $\eps > 0$, is a maximal $\eps$-separated subset
$X$ of $M$ provided with a graph structure given by the sets
$N(p)=\{q \in X \ | \ 0 < d(p,q) < 3 \eps \}$, for any $p \in X$.
In other words, $X$ is such that for any distinct $p$, $q \in X$,
$d(p,q) \geq \eps$ and $\bigcup_{p \in X}B(p,\eps) = M$. Moreover,
$pq$ is an edge if and only if $0 < d(p,q) < 3\eps$. Denote by
$m(p)$ the number of elements of $N(p)$.
\begin{rmk}\label{remark}
    Let us remark that if $B(p,\rho)$ is a
    ball in $M$ with radius $\rho < \frac{1}{2}Inj(M)$, then the
    volume $V(p,\rho)$ of the ball $B(p,\rho)$ is bounded below by
    a constant depending only on $\rho$ and $m$
    (this is Croke's Inequality, see for instance in \cite{Ch}
    p.136). Moreover, if $M$ has Ricci curvature bounded below by
    $-(m-1)\kappa$ then the volume of a ball of radius $R$ is
    bounded above by a constant depending only on $m$, $\kappa$ and $R$
    (this follows from Bishop's comparison theorem, see for
    instance \cite{Ch}, p.126). These bounds will be used frequently in the
sequel.
\end{rmk}
Choose $\eps$ smaller than $\frac{1}{2} Inj (M)$. Denote by
$\kappa \geq 0$ a constant such that $Ricci(M,g) \geq -(m-1)
\kappa g$. Then, using Remark \ref{remark} we can show that $m(p)$
is bounded above by a constant $\nu_X$ depending only on $m$,
$\kappa$ and $\eps $ and that $ \frac{1}{V_{-\kappa}(\eps)} Vol(M)
\leq |X| \leq \frac{2^m}{\eps^m c(m)}Vol(M) $, where
$V_{-\kappa}(\eps)$ denotes the volume of the ball of radius
$\eps$ in the simply connected space of constant sectional
curvature $-\kappa$ and of dimension $m$.

\section{Spectra comparison for rough Laplacian and twisted
Laplacian}\label{comp}

In this section, we will establish the comparison between the
spectra of the rough Laplacian and a twisted Laplacian. Let us
state the main result.

\begin{thm}\label{main}
    Let $m$, $n$ be positive integers,
    $\kappa$, $k_1$, $k_2 \geq 0$ and $r_0 \geq 20 \eps >0$.
    There exist positive constants $c$, $c'$ depending only on $m$, $n$,
    $\kappa$, $k_1$, $k_2$ and $\eps$ such that
    for any $M \in \M$, any vector bundle $E \in \E$ over $M$ satisfying one of the
    following condition
    \begin{enumerate}
        \item[I)] the curvature of $E$ is harmonic i.e. $d^* R^E =
        0$,
        \item[II)] $E$ is of complex (or quaternionic) rank one
    \end{enumerate}
    and for any $\eps$-discretization $X$ of $E$, we can construct a
    canonical twisted Laplacian $\Delta_A$ and a potential $V$ depending only on
    the local geometry of $E$ such that, for $1 \leq k \leq n |X|$
    \begin{eqnarray}
        c' \lambda_k(X,A,V) \leq \lambda_k(E) \leq c \lambda_k(X,A,V) .\nonumber
    \end{eqnarray}
    In particular, if the vector bundle is flat, the potential is
    zero and $\Delta_A$ is a discrete magnetic
    Laplacian.
\end{thm}

Roughly speaking, the basic idea of the proof is the same as to
prove the theorem of comparison of spectra between the Laplacian
acting on functions and the combinatorial Laplacian (\cite{Ma},
Theorem 3.7). But a first fundamental difference between the
functions and the vector bundles cases is the construction of the
twisted Laplacian. Indeed, in \cite{Ma} the combinatorial
Laplacian appearing in Theorem 3.7 is canonically associated to
the graph that discretizes the manifold. For vector bundles, such
a canonical Laplacian on graphs does not obviously exist. Hence, a
first step of the proof consists in constructing a suitable
twisted Laplacian $\Delta_A$ and a potential $V$ (Section
\ref{Constr}) that will depend only on the local geometry. The
construction of $\Delta_A + V$ differs according to the
assumptions $\mathrm{I})$ and $\mathrm{II})$. We will work with
balls centered on $X$ and for both cases the construction of
$\Delta_A$ relies essentially on changes of bases from a ball to a
neighboring ball, but for vector bundles satisfying $\mathrm{II})$
the definition of $\Delta_A$ is slightly harder. A more
significant difference is the construction of the potential $V$.
For rank one vector bundles, $V$ involves only the first
eigenvalue of balls (with Neumann boundary condition), while in
the other case, we will distinguish "small" eigenvalues of balls
from "large" eigenvalues. In rank one vector bundles the $n$ first
eigenvalues (of such a ball) are the same and correspond to the
minimum of the energy, so that it will make easier the estimating
of $V$.

\medskip

After defining the twisted Laplacian and the potential, we follow
the same way of proof as for the case of functions, but the
underlying analysis is much more difficult. For instance, we need
to establish some Sobolev inequalities for sections that requires
fine tools of analysis as Moser's iteration and Sobolev
inequalities for functions (Lemma \ref{lmmApp} in Appendix). The
definition of the smoothing operator $\S$ and the discretizing
operator $\D$ generalizes in some sense the similar operators
defined by Chavel in \cite{Ch} (Sections VI.5.1 and VI.5.2).
Similarly, we establish norms estimations for these operators $\S$
and $\D$ (Propositions \ref{lmmCVBS} and \ref{lmmCVBD}) in order
to compare Rayleigh quotients of sections with Rayleigh quotients
of functions on the discretization. Then, min-max theorem leads to
the result for "small" eigenvalues. It suffices moreover to have
upper bounds on the respective spectra (Lemma \ref{lmmCVBUB}) to
compare "large" eigenvalues and conclude the proof of Theorem
\ref{main} (Section \ref{concl}).

%%%%%%%%%%%%%%%%%%%%%%%%%%%%%%%%%%%%%%%%%%%%%%%%%%%%%%%%%%%%%%%%%%%%%%%%%%%%%%%%%%%%%%%

\subsection{Local extension}\label{Exten}

In this section we define a way to extend a section as parallel as
possible. In the case of flat vector bundles parallel transport is
the suitable tool, because of the lemma below. Let $\tau_{x,p}$
denotes the parallel transport from $E_p$ to $E_x$ along the
minimizing geodesic joining $p$ to $x$ (for $d(p,x) <
\frac{1}{2}Inj(M)$).

\begin{lmm}\label{FlatExt}
    Let $(E,\nabla)$ be a flat Riemannian vector bundle over a
    Riemannian manifold $(M,g)$. Let $p\in M$ and $B(p,r)$ the ball
    centered at $p$ of radius $r < \frac{1}{2} Inj(M)$. Then for any
    $v \in E_p$, the section $\sigma$ over $B(p,r)$ defined by $\sigma(x) =
    \tau_{x,p}v$ is parallel.
\end{lmm}

\textbf{Proof}: see \cite{DK} Section 2.2.1. $\Box$

\medskip

In the non-flat case, extending by parallel transport is not
strong enough for our purpose, because we need to control the
covariant derivative of such extended sections. More precisely, we
want to extend in an energy minimizing way. This means that we
have to take into account local small eigen\-values. Hence, we
introduce eigensections of the Neumann problem on balls which give
an obstruction to extension in a parallel way. Such eigensections
on balls associated to small eigenvalues are almost parallel
(Lemma \ref{LocExten}) and will provide a good way to extend
sections. Nevertheless, it may happen that there are no (or only a
few) small eigensections on a ball. In this case, parallel
transport will be good enough to extend as we will see.
\begin{lmm}\label{LocExten}
    Let $(E, \nabla) \in \E$ over $(M,g) \in
    \M$. For $0 <r < \frac{1}{2}r_0$ and $p \in M$, let
    $\sigma : B(p,r) \to E$ be a section such that
    $\overline{\Delta}\sigma = \lambda \sigma$ for a constant $\lambda \geq 0$.
    Let $0 < \theta <1$.
    Then there exist $0 < c(m) \leq s \leq 1$ and $c$, $c' >0$
    depending on an upper
    bound for $\lambda$ and on $m$, $n$, $\kappa$, $r$, $k_1$, $k_2$ and $\theta$
    such that
    \begin{eqnarray}
        \|\sigma\|_{\infty,\theta r} &\leq& c \|\sigma\|_{2,r} ,\nonumber\\
        \|\nabla\sigma\|_{\infty,\theta r} &\leq& c' \|\nabla\sigma\|^s_{2,r}\nonumber
    \end{eqnarray}
    where $\|\cdot\|_{q,\rho}$ denotes the $L^q$-norm on the ball
    centered at $p$ of radius $\rho$ ($c'$ depends on
    $c\|\sigma\|_{2,r}$ too).

    \medskip

    Moreover, there exists $c''>0$ depending on $c$, $c'$
    and $r$ such that
    \begin{eqnarray}
        |\sigma(x)-\tau_{x,p}\sigma(p)| \leq c'' \|\nabla\sigma\|^s_{2,r}\nonumber
    \end{eqnarray}
    for all $x \in B(p,\theta r)$. If $k_2 = 0$ i.e. if $E$ is of
    harmonic curvature, then $s=1$ in the previous inequalities.
\end{lmm}

\textbf{Proof}: the idea is to use a Moser iteration to prove the
statement. The more technical part of the argument is carried out
in the appendix (see Lemma \ref{lmmApp}). In order to use Lemma
\ref{lmmApp}, let $\delta > 0$ and $u_{\delta} : B(p,r) \to \R$
defined by $u_{\delta} = \sqrt{|\sigma|^2+\delta}$. Then in one
hand $\Delta (u^2_{\delta}) = 2 u_{\delta}\Delta u_{\delta} - 2
|du_{\delta}|^2$ and in the other hand
$\Delta(u_{\delta}^2)=2\langle\sigma,\overline{\Delta}\sigma\rangle
-2 |\nabla \sigma|^2$ which implies that
\begin{eqnarray}
    u_{\delta}\Delta u_{\delta} \leq \langle
    \sigma,\overline{\Delta}\sigma\ \rangle = \lambda|\sigma|^2 \leq
    \lambda u_{\delta}^2 .\nonumber
\end{eqnarray}
We can then apply Lemma \ref{lmmApp} to $u_{\delta}$ and we get
that $\|u_{\delta}\|_{\infty,\theta r} \leq c
\|u_{\delta}\|_{2,r}\nonumber$. Then let $\delta \to 0$ to obtain
the first claim.

\medskip

For the second inequality, let $\delta > 0$ and $v_{\delta} :
B(p,r) \to \R$ defined by $v_{\delta}(x) = \sqrt{|\nabla
\sigma(x)|^2+\delta}$. Then $$\Delta(v_{\delta}^2) =
    2v_{\delta}\Delta v_{\delta}- 2|dv_{\delta}|^2 =
    2\langle\nabla\sigma,\overline{\Delta}(\nabla \sigma) \rangle
    - 2 |\nabla \nabla \sigma|^2 .$$
But we have that $|\nabla\nabla \sigma|^2-|dv_{\delta}|^2 \geq 0$
and therefore
\begin{eqnarray}
    v_{\delta}\Delta v_{\delta} \leq \langle
    \nabla\sigma,\overline{\Delta}(\nabla \sigma) \rangle =
    \langle\nabla\sigma,\overline{\Delta}(\nabla
    \sigma)-\nabla(\overline{\Delta}\sigma) \rangle +
    \lambda|\nabla \sigma|^2 .
    \nonumber
\end{eqnarray}
By a commuting argument (see \cite{ACGR}, Lemma 2.3) we have for a
local orthonormal frame $\{X_i\}_{i=1,\ldots,m}$ of $M$
\begin{multline}
    \langle\nabla\sigma,\overline{\Delta}(\nabla\sigma)-\nabla(\overline{\Delta}\sigma)\rangle =\\
    \lambda|\nabla \sigma|^2 -
        \langle \nabla_{Ric(\cdot)}\sigma,\nabla\sigma \rangle - 2
        \sum_{i=1}^m \langle R^E(X_i,\cdot) \nabla_{X_i}\sigma,\nabla \sigma\rangle
        + \langle (d^* R^E)\sigma , \nabla \sigma \rangle
    \nonumber
\end{multline}
and as $Ricci(M,g) \geq -(m-1)\kappa g$, $|R^E|\leq k_1$ and $|d^*
R^E| \leq k_2$ we then get
\begin{eqnarray}
    \langle\nabla\sigma,\overline{\Delta}(\nabla\sigma)-\nabla(\overline{\Delta}\sigma)\rangle
    \leq \left(\lambda + (m-1)\kappa + 2n^2 k_1\right)|\nabla
    \sigma|^2 + n^2 k_2 |\sigma| |\nabla \sigma| .\nonumber
\end{eqnarray}
By the first part of the proof, we obtain that on $B(p,\theta r)$
\begin{multline}
    \langle\nabla\sigma,\overline{\Delta}(\nabla\sigma)-\nabla(\overline{\Delta}\sigma)\rangle
    \leq \\ \left(\lambda + (m-1)\kappa + 2n^2 k_1\right)|\nabla
    \sigma|^2 + n^2 k_2 c \|\sigma\|_{2,r} |\nabla \sigma| \nonumber
\end{multline}
and this implies (on $B(p,\theta r))$
\begin{eqnarray}
    \Delta v_{\delta} \leq \left(\lambda + (m-1)\kappa +
    2n^2 k_1\right) v_{\delta} + n^2 k_2 c \|\sigma\|_{2,r} .
    \nonumber
\end{eqnarray}
If $\theta' < \theta$ we can apply Lemma \ref{lmmApp} to
$v_{\delta}$ and let $\delta \to 0$ to obtain
\begin{eqnarray}
    \|\nabla \sigma \|_{\infty,\theta' r}
    \leq
    c' \|\nabla \sigma\|^s_{2,\theta r}
     \leq
    c' \|\nabla \sigma\|^s_{2,r} .\label{deux}
\end{eqnarray}
Note that if $k_2=0$, then $s=1$ and $c'$ does not depend on $c
\|\sigma\|_{2,r}$. The two first inequalities in the statement are
then true for any $\theta'$ such that $0 < \theta' < \theta < 1$.
So rename $\theta'$ by $\theta$ to obtain the statement.

\medskip

Finally, recall that if $\gamma$ is the minimizing geodesic
joining $p$ to $x \in B(p,\theta r) $ of length $l$ ($ < \theta
r$), then $|\sigma(x)-\tau_{x,p}\sigma(p)| \leq \int_0^l
\left|\nabla_{\dot{\gamma}(t)} \sigma(\gamma(t))\right|dt \leq l
\|\nabla \sigma\|_{\infty,\theta r} $. Using (\ref{deux}) leads to
the result. $\Box$

%%%%%%%%%%%%%%%%%%%%%%%%%%%%%%%%%%%%%%%%%%%%%%%%%%%%%%%%%%%%%%%%%%%%%%%%%%%%%%%%%%%%%%%

\medskip

From now on, let $E \in \E$ over $M \in \M$ and fix $\eps \leq
\frac{1}{20}r_0$. Let $X$ be an $\eps$-discretization of $E$. Let
$\sigma_k^p : B(p,10\eps) \to E$ be the eigensection associated to
the $k^{th}$ eigenvalue $\lambda_k(p)$ of $\overline{\Delta}$ on
$B(p,10\eps)$ with Neumann boundary condition such that
$\int_{B(p,10\eps)}\langle\sigma_k^p,\sigma_l^p\rangle dV =
\delta_{kl}V(p,10\eps)$.

\begin{rmk}
    If $E$ is flat $\lambda_1(p) = \ldots = \lambda_n(p) =0$ and the
    $\sigma_k^p$'s give a local orthonormal frame over $B(p,10
    \eps)$.
\end{rmk}

\begin{rmk}
    If $n=2$ (resp. $n=4$) and $E$ is of complex (resp. quaternionic) rank
    one, then $\lambda_1(p)=\ldots = \lambda_n(p)$. Indeed, the
    section $i \sigma_1^p$ (resp. $i \sigma_1^p$, $j \sigma_1^p$,
    $k \sigma_1^p$ where $i$, $j$, $k$ are the quaternions with
    $i^2 = j^2 = k^2 =-1$) satisfies $\nabla(i \sigma_1^p) = i
    \nabla \sigma_1^p$ which implies that $i \sigma_1^p$ is a
    $\lambda_1(p)$-eigensection orthogonal to $\sigma_1^p$. Hence,
    we can choose $\sigma_k^p$ such that for any $x$ in $B(p,10\eps)$,\linebreak
    $\langle \sigma_k^p(x),\sigma_l^p(x) \rangle =0$ for any $1 \leq k\leq n$,
    $1 \leq l \leq n$, $k\neq l$.
\end{rmk}

\begin{lmm}\label{lmmDel}
    Let $0 \leq \alpha < \frac{1}{n+1}$. There exists $\delta >0$
    depending only on $\alpha$, $m$, $n$, $k_1$, $k_2$, $\kappa$,
    $\eps$ such that if $\lambda_k(p) \leq \delta$ then $\forall \ 1\leq i$,
    $j \leq k$ and $\forall x \in B(p,8\eps)$
    $|\langle\sigma_i^p(x),\sigma_j^p(x) \rangle - \delta_{ij}| \leq
    \alpha $. In particular, if $\lambda_k(p) \leq \delta$, then
    $\{\sigma_1^p(x), \ldots, \sigma_k^p(x)\}$ spans a
    $k$-dimensional vector subspace of $E_x$, for any $x\in B(p,8\eps)$.
\end{lmm}

To prove this lemma, let us recall a basic fact of linear algebra
(the proof of the fact is left to the reader). Let $V$ be an
$n$-dimensional vector space provided with an inner product
$\langle \cdot , \cdot \rangle$. If $\{v_1, \ldots, v_n \}
\subseteq V$ is such that $\left| \langle v_i, v_j \rangle -
\delta_{ij} \right| \leq \alpha < \frac{1}{n+1}$, then $\{v_1,
\ldots, v_n \}$ is a basis of $V$. Moreover for any
$v=\sum_{i=1}^n a_i v_i$, we have $(1 - \alpha(n+1)) \sum_{i=1}^n
a_i^2 \leq \|v\|^2 \leq (1 + \alpha(n+1)) \sum_{i=1}^n a_i^2$.
Such a basis will be referred as an \textbf{almost orthonormal
basis}.

\medskip

\textbf{Proof of Lemma \ref{lmmDel}}: let $ f_{ij}(x) = \langle
\sigma^p_i(x),\sigma^p_j(x) \rangle $ and denote by $m_{ij}$ its
mean over $B(p,10\eps)$, then
\begin{eqnarray}
    m_{ij} = \frac{1}{V(p,10\eps)} \int_{B(p,10\eps)} f_{ij} dV =
    \delta_{ij} .\nonumber
\end{eqnarray}
A result of Kanai ensuring the existence of $c_K >0$ depending
only on $\eps$ and $\kappa$ (see \cite{Ch}, Lemma VI.5.5) and the
assumption $\lambda_k(p) \leq \delta$ imply
\begin{eqnarray}
    0 \leq \int_{B(p,10\eps)} |f_{ij} - \delta_{ij}|dV  \leq   c_K
    \int_{B(p,10\eps)}|df_{ij}| dV
      \leq  c_K V(p,10 \eps) \sqrt{\delta} .\label{eq61}
\end{eqnarray}
Moreover,
\begin{eqnarray}
    \inf_{x \in B(p,\frac{\eps}{2})} \{|f_{ij}(x) - \delta_{ij}|\} V\left(p,\frac{\eps}{2}\right)
    & \leq & \int_{B(p,\frac{\eps}{2})} |f_{ij}(x) - \delta_{ij}|dV(x)
    \nonumber \\
    & \leq & c_K V(p,10 \eps) \sqrt{\delta} .\label{eq62}
\end{eqnarray}

The last inequality follows from (\ref{eq61}). Hence (\ref{eq62})
implies that there exists $p' \in M$, $d(p,p') \leq
\frac{\eps}{2}$, such that
\begin{eqnarray}
    |\langle \sigma^p_{i}(p'),\sigma^p_j(p')\rangle - \delta_{ij}|
    \leq 2 c_K \frac{V(p,10\eps)}{V(p,\frac{\eps}{2})}\sqrt{\delta}
    \leq c \sqrt{\delta} .\nonumber
\end{eqnarray}
We conclude then as follows
\begin{multline}
    \left|\langle \sigma^p_{i}(x),\sigma^p_j(x)\rangle - \delta_{ij}\right|
    \leq \\
    \left|\langle \sigma^p_{i}(x),\sigma^p_j(x)\rangle -
    \langle \tau_{x,p'}\sigma^p_{i}(p'),\tau_{x,p'}\sigma^p_j(p')
    \rangle\right| +
    \left|\langle \sigma^p_{i}(p'),\sigma^p_j(p')\rangle-\delta_{ij}\right| \\
    \leq \left|\langle \sigma^p_{i}(x),\sigma^p_j(x)\rangle - \langle
    \tau_{x,p'}\sigma^p_{i}(p'),\tau_{x,p'}\sigma^p_j(p')\rangle\right| +
    c \sqrt{\delta} .\label{eq63}
\end{multline}
For any $x \in B(p, 8 \eps)$ the minimizing geodesic
$\overline{xp'}$ stays in $B(p,9 \eps)$, so we can write
\begin{multline}
    \left|\langle \sigma^p_{i}(x),\sigma^p_j(x)\rangle - \langle\sigma^p_{i}(p'),\sigma^p_j(p')
    \rangle\right| \leq 9 \eps \|d\langle \sigma^p_{i},\sigma^p_j \rangle\|_{\infty,9 \eps}\\
    \leq 9 \eps \left(\|\nabla \sigma^p_{i}\|_{\infty,9 \eps} \|\sigma^p_j\|_{\infty,9\eps}
    + \| \sigma^p_{i}\|_{\infty,9 \eps} \|\nabla \sigma^p_j \|_{\infty,9 \eps}\right)\\
    \leq 9 \eps c' \left(\|\nabla \sigma^p_{i}\|^s_{2,10\eps}\|\sigma^p_j\|_{2,10\eps}
    + \| \sigma^p_{i}\|_{2,10\eps} \|\nabla \sigma^p_j \|^s_{2,10\eps}\right)
    \nonumber
\end{multline}
where we used Lemma \ref{LocExten} in the last inequality. By
definition of the $\sigma^p_i$'s and by assumption on
$\lambda_i(p)$ we get
\begin{eqnarray}
    |\langle \sigma^p_{i}(x),\sigma^p_j(x)\rangle - \langle\sigma^p_{i}(p'),\sigma^p_j(p') \rangle|
    \leq c'' \sqrt{\delta^s} .\label{eq64}
\end{eqnarray}
Finally, (\ref{eq63}) and (\ref{eq64}) imply that for a
sufficiently small $\delta$ we have
\begin{eqnarray}
    |\langle \sigma^p_{i}(x),\sigma^p_j(x)\rangle - \delta_{ij}|
    \leq \left(c \sqrt{\delta} + c''\sqrt{\delta^{s}}\right)
    \leq \alpha < \frac{1}{n+1} \nonumber
\end{eqnarray} and this ends the proof. $\Box$
\begin{dfn}
    Fix once and for all $0< \alpha < \frac{1}{n+1}$.
    Let $\delta$ be given by Lemma \ref{lmmDel}. For $p \in X$, define then
    {\boldmath $\mu(p)$} as the largest integer such that $\lambda_{\mu(p)}(p) \leq \delta$.
\end{dfn}
\begin{rmk}
    If the vector bundle is flat, $\mu(p) = n$, for any $ p \in X$.
\end{rmk}
For $p \in X$, we want to extend a section in a neighborhood of
$p$ as parallel as possible and taking care of local small
eigenvalues as said before. So let us define the \textbf{local
extension} that associates to a vector in $E_p$ a local section
over $B(p,10 \eps)$. Consider $E_{\mu(p)}$ the
$\mu(p)$-dimensional vector subspace of $E_p$ spanned by
$\{\sigma_1^p(p),\ldots,\sigma_{\mu(p)}^p(p)\}$. Let
$E_{\mu(p)}^{\bot}$ be the orthogonal complement of $E_{\mu(p)}$
in $E_p$ and choose $\{e_{\mu(p)+1}^p,\ldots,e_n^p\}$ an
orthonormal basis of $E_{\mu(p)}^{\bot}$. By construction,
$\{e_1^p = \sigma_1^p(p),\ldots,e_{\mu(p)}^p =
\sigma_{\mu(p)}^p(p),e_{\mu(p)+1}^p,\ldots,e_n^p\}$ is an almost
orthonormal basis of $E_p$. We extend this basis on $B(p,10\eps)$
by
\begin{eqnarray}
     e_i^p(x) :=
     \left\{    \begin{array}{ll}
                    \sigma_i^p(x) & \text{if } i\leq \mu(p) ,\\
                    \tau_{x,p}e_i^p & \text{otherwise}
                \end{array}\right.        \nonumber
\end{eqnarray}
and we define the local extension of $v = \sum_{i=1}^n v_i e_i^p$
by $\sum_{i=1}^n v_i e_i^p(x)$.
\begin{rmk}
    If $E$ is flat, the local extension corresponds to the extension by parallel
    transport along radial geodesics. In this case, it suffices to
    choose any orthonormal basis $\{e_1^p, \ldots, e_n^p\}$ of
    $E_p$ and extend it radially to obtain $\{e_1^p(x), \ldots,
    e_n^p(x)\}$.
\end{rmk}
\begin{lmm}
    For any $x \in B(p,8\eps)$, $\{e_1^p(x),\ldots,e_n^p(x)\}$ is
    an almost orthonormal basis of $E_x$.
\end{lmm}

\textbf{Proof}: if $\mu(p)=0$ the claim is clearly true. If
$\mu(p)=n$ the claim follows from Lemma \ref{lmmDel}. Hence
suppose $1 \leq \mu(p) \leq n-1$. By Lemma \ref{lmmDel} $\langle
e_1^p(x), \ldots,e_{\mu(p)}^p(x) \rangle$ is $\mu(p)$-dimensional
and  as parallel translation preserves the inner product $\langle
e_{\mu(p)+1}^p(x),\ldots,e_n^p(x) \rangle$ is
$(n-\mu(p))$-dimensional. So we have to show that there exists $c
>0$ such that
\begin{eqnarray}
    |\langle e_i^p(x),e_j^p(x) \rangle |\leq c < \frac{1}{n+1} \ ,\ \ \
    \ \forall 1\leq i \leq \mu(p) < j \leq n.
    \nonumber
\end{eqnarray}
Let us prove this estimate. As $e_j^p(p)$ and $\sigma_i^p(p)$ are
orthogonal, we have
\begin{eqnarray}
    \left|\langle e_j^p(x),e_i^p(x) \rangle\right| & = &
    \langle e_j^p(p),\tau_{p,x}\sigma_i^p(x) - \sigma_i^p(p) \rangle\nonumber\\
    & \leq &
    \left|e_j^p(p)\right|\cdot\left|\sigma_i^p(x)- \tau_{x,p}\sigma_i^p(p)\right|\nonumber
    = \left|\sigma_i^p(x)- \tau_{x,p}\sigma_i^p(p)\right| .
\end{eqnarray}
By Lemma \ref{LocExten} $|\sigma_i^p(x)- \tau_{x,p}\sigma_i^p(p)|
\leq c' \sqrt{\delta^s}$. Hence $\left|\langle e_j^p(x),e_i^p(x)
\rangle\right| \leq c' \sqrt{\delta^s}$. Then, readjust $\delta$
if necessary to obtain $\left|\langle e_j^p(x),e_i^p(x)
\rangle\right| \leq c < \frac{1}{n+1}$. $\Box$
\begin{rmk}
    For the sequel, let $\delta'$ denote a constant, $0< \delta' <1$, such
    that $(1-\delta')\sum_{i=1}^n v_i^2 \leq \left| \sum_{i=1}^n v_i
    e_i^p(x)\right|^2 \leq (1+\delta')\sum_{i=1}^n v_i^2$, for any $x \in B(p,8\eps)$.
\end{rmk}
\begin{lmm}\label{lmmDev}
    There exists a positive constant $c$ depending only on $n$, $k_1$, $\eps$
    such that for any $p\in X$ and any $\mu(p) < i  \leq n$,
    $\|\nabla e_i^p\|_{\infty,9\eps}\leq c$.
\end{lmm}

\textbf{Proof}: let $x \in B(p,9\eps)$ and consider $\gamma$ the
minimizing geodesic from $p$ to $x$ of length $l$ ($l <9\eps$) and
$\{X_1 = \dot{\gamma}(t), X_2, \ldots , X_n\}$ an orthonormal
basis of $E_x$ with $\nabla_{X_i}X_j=0$. Then
\begin{eqnarray}
    \left|\nabla e_i^p(x)\right| ^2 = \sum_{j=1}^n \left|\nabla_{X_j}e_i^p(x)\right|^2 \leq
    \sum_{j=1}^n\left(
    \int_{0}^{l}\left|\nabla_{\dot{\gamma}(t)}\nabla_{X_j}e_i^p(x)\right|
    dt\right)^2
    \nonumber
\end{eqnarray}
but $|R^E(\dot{\gamma}(t),X_j)e_i^p|
=|\nabla_{\dot{\gamma}(t)}\nabla_{X_j}e_i^p| \leq k_1 $. Therefore
$|\nabla e_i^p(x)|^2 \leq k_1^2 l^2 n$ and this concludes the
proof. $\Box$

%%%%%%%%%%%%%%%%%%%%%%%%%%%%%%%%%%%%%%%%%%%%%%%%%%%%%%%%%%%%%%%%%%%%%%%%%%%%%%%%%%%%%%%

\subsection{Construction of the twisted Laplacian} \label{Constr}

The construction of $\Delta_A$ differs according to the
assumptions done on $E$. However, the basic idea is the same in
all cases and relies on the fact that $A$ has to express the
holonomy. So let us consider $p$,$q \in X$, $p \in N(q)$ and let
$x \in B(p,8\eps)\cap B(q,8\eps)$. Then define $a(p,q)_{ij}(x)$ by
\begin{eqnarray}
    e_j^p(x)=\sum_{i=1}^n a(p,q)_{ij}(x)e_i^q(x) \ \ \ \ \ \ \ \forall
    j=1, \ldots,n
    \nonumber
\end{eqnarray}
where $e_i^p$, $e_j^q$ are defined in Section \ref{Exten}. We
define $A(p,q) : \R^n \to \R^n$ on the canonical basis
$\{e_1,\ldots e_n\}$ of $\R^n$ by $A(p,q)e_j = \sum_{i=1}^n
A(p,q)_{ij} e_i$, where $A(p,q)_{ij}$ is defined as follows.

\medskip

\textbf{If $\mathbf{E}$ is of harmonic curvature} then define
$A(p,q)_{ij}$ by
\begin{eqnarray}
    A(p,q)_{ij} = a(p,q)_{ij}(q)\nonumber
\end{eqnarray}

\textbf{If $\mathbf{E}$ is of complex (or quaternionic) rank one}
then define $A(p,q)_{ij}$ by
\begin{eqnarray}
    A(p,q)_{ij} = \frac{1}{V_{pq}}\int_{B_{pq}}a(p,q)_{ij}(x)dV(x)\nonumber
\end{eqnarray}
where $B_{pq}$ is the ball centered at the mid-point of $p$ and
$q$ of radius $5 \eps$ and $V_{pq}$ denotes its volume. Note that
$B_{pq} \supseteq B(p,3\eps) \cup B(q,3\eps)$.
\begin{rmk}\label{rmk321}
In the canonical basis of $\R^n$, we can write
\begin{eqnarray}
    D_A f(p,q) = \sum_{i=1}^n D_A f(p,q)_i e_i = \sum_{i=1}^n
    \left(f_i(q)-\sum_{j=1}^n A(p,q)_{ij}f_j(p) \right) e_i
    \nonumber
\end{eqnarray}
\end{rmk}
\begin{rmk}\label{rmk432}
If $E$ is of harmonic curvature, we have by definition \\
$e_j^p(q)=\sum_{i=1}^n A(p,q)_{ij}e_i^q(q)$, $\forall j=1,
\ldots,n$.
\end{rmk}
\begin{rmk}
If $E$ is flat, $a(p,q)_{ij}(x)$ is constant and so for $ j=1,
\ldots,n$ and for any $x \in B(p,8\eps) \cap B(q,8\eps)$,
$e_j^p(x)=\sum_{i=1}^n A(p,q)_{ij}e_i^q(x)$. Moreover, in this
case $A(p,q)A(p,q)^t = Id$ and $A(p,q)^t=A(q,p)$. So that
$\Delta_A$ is a discrete magnetic Laplacian.
\end{rmk}
\textbf{If $\mathbf{E}$ is of harmonic curvature} let $V :\FnX \to
\FnX$ be defined by
\begin{eqnarray}
        (Vf)(p) = \sum_{i\leq \mu(p)} \lambda_i(p) f_i(p)e_i +
        \sum_{i>\mu(p)}f_i(p)e_i .
        \nonumber
\end{eqnarray}
\textbf{If $\mathbf{E}$ is of complex (or quaternionic) rank one}
let $V :\FnX \to\FnX$ be defined by
\begin{eqnarray}
    (Vf)(p) = \left( \lambda_1(p) + \sum_{q \in
    N(p)}\lambda_1(q)\right)f(p).
\end{eqnarray}
\begin{rmk}
    If the vector bundle is flat, then we have $V=0$.
\end{rmk}
%%%%%%%%%%%%%%%%%%%%%%%%%%%%%%%%%%%%%%%%%%%%%%%%%%%%%%%%%%%%%%%%%%%%%%%%%%%%%%%%%%%%%%%
\subsection{Smoothing operator}

\begin{dfn}
    Let $\{\psi_p\}_{p\in X}$ be a partition of unity subordinate to
    the cover $\{B(p,2\eps)\}_{p\in X}$. Define the \textbf{smoothing
    operator}
    $\S : \FnX \to \Gamma(E)$ by
    \begin{eqnarray}
        (\S f)(x) = \sum_{p\in X} \psi_p(x) \left( \sum_{i=1}^n f_i(p)
        e_i^p(x) \right)
        \nonumber
    \end{eqnarray}
    where $f(p) = \sum_{i=1}^n f_i(p)e_i$.
\end{dfn}
\begin{prp}\label{lmmCVBS}
    There exist constants $c_0$, $c_1$, $c_2$ and $\Lambda > 0$
    depending only on $m$, $n$, $k_1$, $k_2$, $\kappa$ and $\eps$ such that
    \begin{enumerate}
        \item[i)] $\forall f \in \FnX$, $\|\S f\|^2 \leq c_0
        \|f\|^2$,
        \item[ii)]$\forall f \in \FnX$, $\|\nabla(\S f)\|^2 \leq c_1\left( \|D_A f\|^2 +
        (Vf,f)\right)$,
        \item[iii)]$\forall f \in \FnX$ with $\|D_A f\|^2 + (Vf,f) \leq \Lambda
        \|f\|^2$, $\|\S f\|^2 \geq c_2 \|f\|^2$ holds.
    \end{enumerate}
\end{prp}
\textbf{Proof}: for the first inequality note that
$\{B(p,\eps)\}_{p\in X}$ covers $M$. Hence
\begin{eqnarray}
    \|\S f\|^2 & \leq &
    \sum_{q \in X} \int\limits_{B(q,\eps)}
    \left|\sum_{p\in B(q,3\eps)\cap X}\psi_p(x)\sum_{i=1}^n f_i(p)e_i^p(x)\right|^2 dV(x)\nonumber\\
    & \leq &
     (1+\delta') \sum_{q \in X}V(q,\eps)\!\!\!
    \sum_{p\in B(q,3\eps)\cap X}\!\! |f(p)|^2
    \ \leq\  (1+\delta') c \|f\|^2 .\nonumber
\end{eqnarray}
In order to prove $ii)$ fix $q \in X$ and let $x \in B(q,\eps)$.
Then as $\{\psi_p\}_{p \in X}$ is a partition of unity, we have
$\sum_{p \in X} d\psi_p = 0$, so that we can write
\begin{eqnarray}
    \nabla (\S f)(x) & = &
    \sum_{p\in B(q,3\eps)\cap X} \psi_p(x)\Big(\sum_{i=1}^n f_i(p)\nabla
    e_i^p(x)\Big) +  \nonumber\\
    & & \sum_{p\in N(q)}d \psi_p(x)
    \Big(\sum_{i=1}^n f_i(p)e_i^p(x) - \sum_{i=1}^n f_i(q)e_i^q(x)
    \Big) .\label{nabla4}
\end{eqnarray}
Then, Lemma \ref{lmmDev} implies
\begin{multline}
    \int_{B(q,\eps)}\left|
    \sum_{p\in B(q,3\eps)\cap X} \psi_p(x)\Big(\sum_{i=1}^n f_i(p)\nabla
    e_i^p(x)\Big)\right|^2dV(x) \leq \\
    n \sum_{p\in B(q,3\eps)\cap X}
    \left(\sum_{i\leq \mu(p)} f_i(p)^2 \int_{B(q,\eps)}|\nabla
           e_i^p(x)|^2 dV(x) + c \sum_{i > \mu(p)} f_i(p)^2
           \right)\\
       \leq c' \sum_{p\in B(q,3\eps)\cap X}\left(Vf\right)(p)\cdot f(p) .\label{3star}
\end{multline}

To estimate the second term of (\ref{nabla4}), we need the
following lemma.
\begin{lmm}\label{lmm320}
There exists a positive constant $c$ depending only on $m$, $n$,
$k_1$, $k_2$, $\kappa$ and $\eps$  such that
\begin{multline}
    \int_{B(q,\eps)}\left|\sum_{i=1}^n f_i(p)e_i^p(x) - \sum_{i=1}^n
    f_i(q)e_i^q(x)\right|^2 \leq
    \\ c \left(|D_Af(q,p)|^2 + (Vf)(p)\cdot f(p) + (Vf)(q)\cdot
    f(q)\right) .\nonumber
\end{multline}
\end{lmm}

\textbf{Proof}: see Appendix \ref{applmm1}. $\Box$

\medskip

Hence by (\ref{3star}), (\ref{nabla4}) and Lemma \ref{lmm320} we
get
\begin{multline}
    \int_{B(q,\eps)}|\nabla(\S f)(x)|^2dV(x) \leq \\
    c'' \sum_{p\in
    B(q,3\eps)\cap X}\Big(|D_A f(q,p)|^2 + (Vf)(p)\cdot f(p)
     + (Vf)(q)\cdot f(q)\Big) .
    \nonumber
\end{multline}
Then summing on $q \in X$ implies the claim.

\medskip

To prove the third part of Proposition \ref{lmmCVBS}, define
$(\S_q f)(x) = \sum\limits_{i=1}^n f_i(q)e_i^q(x)$ for $x$ in
$B(q,\frac{\eps}{2})$. Then, by Lemma \ref{lmm320} we get
\begin{multline}
    \int\limits_{B(q,\frac{\eps}{2})}\left|(\S f )(x) - (\S_q f)(x)
    \right|^2dV(x) = \\
    \int\limits_{B(q,\frac{\eps}{2})}\left| \sum_{\ p\in N(q)}
    \psi_p(x)\sum_{j=1}^n\big(f_j(p)e_j^p(x) - f_j(q)e_j^q(x) \big)\right|^2 dV(x) \leq \\
    c \sum_{p \in N(q)} \left (\left| D_A
    f(q,p)\right|^2 + (Vf)(p)\cdot f(p) + (Vf)(q)\cdot
    f(q)\right) .
    \label{eqn2}
\end{multline}
As the balls of radius $\frac{\eps}{2}$ centered on $X$ are
disjoint, we can write
\begin{eqnarray}
    \| \S f \|^2 %& \geq &
    %\sum_{q \in X} \int\limits_{B(q,\frac{\eps}{2})}\left| \S f(x)\right|^2dV(x)\nonumber\\
    & \geq &
    \sum_{q \in X} \int\limits_{B(q,\frac{\eps}{2})}\left|
    \left(\S_q f(x) -\S f(x)\right)  - \S_q f(x)
    \right|^2dV(x)\nonumber \\
    & \geq &
    \sum_{q \in X} \int\limits_{B(q,\frac{\eps}{2})}\left| \S_q
    f(x)\right|^2dV(x)\nonumber \\
    &  &
    - 2 \sum_{q \in X} \int\limits_{B(q,\frac{\eps}{2})}
    \left| \S_q f(x)\right| \left|\S f(x) -\S_q f(x) \right|dV(x) .\nonumber
\end{eqnarray}
By construction, $(1-\delta') |f(q)|^2\leq \left| \S_q
f(x)\right|^2 \leq (1+\delta')|f(q)|^2$ and by Cauchy-Schwarz
inequality combined with (\ref{eqn2}), we get
\begin{multline}
    \sum_{q \in X} \int\limits_{B(q,\frac{\eps}{2})}
    \left| \S_q f(x)\right| \left|\S f(x) -\S_q f(x) \right|dV(x) \\ \leq
    c'(1+\delta') \|f\|\sqrt{ \|D_A f\|^2+(Vf,f)} .\nonumber
\end{multline}
Hence, $\|\S f\|^2 \geq (1-\delta')c'' \|f\|^2 -
2c'(1+\delta')\|f\|\sqrt{ \|D_A f\|^2+(Vf,f)}$. Choose $\Lambda
> 0$ sufficiently small so that if $f$ satisfies $\|D_A f\|^2 +
\left( Vf , f \right) \leq \Lambda \|f\|^2$, then
\begin{eqnarray}
    \| \S f \|^2  \geq
    \|f\|^2  \left((1-\delta')c''-2c'(1+\delta')\sqrt{\Lambda} \right) \geq
    \frac{(1-\delta')c''}{2}\|f\|^2 .\nonumber
\end{eqnarray}
This concludes the proof of Proposition \ref{lmmCVBS}.  $\Box$

%%%%%%%%%%%%%%%%%%%%%%%%%%%%%%%%%%%%%%%%%%%%%%%%%%%%%%%%%%%%%%%%%%%%%%%%%%%%%%%%%%

\subsection{Discretizing operator}

\begin{dfn}
    Define the \textbf{discretizing operator} $\D : \Gamma(E) \to \FnX$
    by
    \begin{eqnarray}
        (\D s)(p) = \sum_{i=1}^n \frac{1}{V(p,3\eps)}
        \int_{B(p,3\eps)} s_i^p(x) dV(x) e_i
        \nonumber
    \end{eqnarray}
    where $s(x) = \sum_{i=1}^n s_i^p(x)e_i^p(x)$ for $x$ in $B(p,3\eps)$.
\end{dfn}
\begin{prp}\label{lmmCVBD}
    There exist constants $c'_0$, $c'_1$, $c'_2$ and $\Lambda' > 0$
    depending only on $m$, $n$, $\kappa$, $k_1$, $k_2$ and $\eps$ such that
    \begin{enumerate}
        \item[i)] $\forall s \in \Gamma(E)$, $\|\D s\|^2 \leq c'_0 \|s\|^2$,
        \item[ii)]$\forall s \in \Gamma(E)$, $\|D_A (\D s)\|^2 +
        (V(\D s),\D s) \leq c'_1 \|\nabla s\|^2$,
        \item[iii)]$\forall s \in \Gamma(E)$ such that $\|\nabla s\|^2 \leq
        \Lambda' \|s\|^2$, $\|\D s\|^2 \geq c'_2 \|s\|^2$ holds.
    \end{enumerate}
\end{prp}

\textbf{Proof}: the first point follows directly from the
following inequality
\begin{eqnarray}
    |\D s(p)|^2 \leq c
    \int_{B(p,3\eps)} \sum_{i=1}^n |s_i^p(x)|^2 dV(x)
    \leq c(1-\delta')^{-1} \int_{B(p,3\eps)}|s(x)|^2dV(x) .\nonumber
\end{eqnarray}
To prove the second point, we first prove that
\begin{eqnarray}
    \|D_A(\D s)\|^2 + (V(\D s),\D s) & \leq &
    c \left(\|\nabla s\|^2 +
    \sum_{p\in X}\left(\widetilde{V}s\right)(p)\right) \label{eqn13}
\end{eqnarray}
where if \textbf{$\mathbf{E}$ is of harmonic curvature} then
\begin{eqnarray}
    \left(\widetilde{V}s\right)(p) =
    \left( \sum_{i\leq \mu(p)}\lambda_i(p)\int_{B(p,3\eps)}|s_i^p|^2dV +
    \sum_{i > \mu(p)} \int_{B(p,3\eps)}|s_i^p|^2dV \right)
    \nonumber
\end{eqnarray}
and if \textbf{$\mathbf{E}$ is of complex (or quaternionic) rank
one}
\begin{eqnarray}
    \left(\widetilde{V}s\right)(p) =
    \left(\lambda_1(p)+ \sum_{q \in N(p)}\lambda_1(q)\right) \int_{B(p,3\eps)}|s|^2dV
    \nonumber
\end{eqnarray}
and $s$ is written locally as $s(x) = \sum_{i=1}^n
s_i^p(x)e_i^p(x)$ for $x \in B(p,8\eps)$. First, $\left|(\D
s)(p)_j\right|^2 \leq c \int_{B(p,3\eps)}|s_j^p(x)|^2dV(x)$
implies obviously
\begin{eqnarray}
    (V(\D s),\D s) &\leq& \sum_{p\in X} c'\left(\widetilde{V}s\right)(p) .\label{eqn15}
\end{eqnarray}
Secondly, for $p$ and $q \in N(p)$ let us introduce $B'_{pq}
\subseteq B(p,3\eps)\cap B(q,3\eps)$ the ball centered at the
mid-point of $p$ and $q$ of radius $\eps$ and $V'_{pq}$ its
volume. Then
\begin{eqnarray}
   \lefteqn{|D_A (\D s)(q,p)|^2  =} \nonumber\\
    &&\sum_{i=1}^n \left( \frac{1}{V'_{pq}} \int_{B'_{pq}}
    \left| \D s(p)_i - \sum_{j=1}^n A(q,p)_{ij}\D s(q)_j
    \right| dV(y)
    \right)^2
    \nonumber\\
    &\leq&
    3\sum_{i=1}^n \left( \frac{1}{V'_{pq}}\int_{B'_{pq}}\left|
        \D s(p)_i - s_i^p(y)
    \right|dV(y)\right)^2 \label{eqn8}\\
    &+&
    3\sum_{i=1}^n \left( \frac{1}{V'_{pq}}\int_{B'_{pq}}\left|
        \sum_{j=1}^n A(q,p)_{ij} \left(s_j^q(y)-\D s(q)_j \right)
    \right|dV(y)\right)^2 \label{eqn7}\\
    &+&
    3\sum_{i=1}^n \left( \frac{1}{V'_{pq}}\int_{B'_{pq}}\left|
        s_i^p(y)-\sum_{j=1}^n A(q,p)_{ij} s_j^q(y)
    \right|dV(y)\right)^2 .\label{eqn9}
\end{eqnarray}
We estimate each of the three terms separately.

\medskip

By a result of Kanai (see \cite{Ch}, Lemma VI.5.5), there exists
$c_K >0$ depending only on $\eps$ and $\kappa$ such that
\begin{eqnarray}
    \frac{1}{V'_{pq}}\int_{B'_{pq}} \left| \D s(p)_i - s_i^p(y) \right|dV(y)
    &\leq&
    c_K \int_{B(p,3 \eps)} \left|d s_i^p(y)\right|dV(y).\nonumber
\end{eqnarray}
Moreover
\begin{eqnarray}
     \sqrt{1-\delta'} \left|d s_i^p(y)\right| \leq \left|\sum_{j=1}^n d s_j^p(y)e_j^p(y)\right| =
    \left|\nabla s(y) - \sum_{j=1}^n s_j^p(y)\nabla e_j^p(y)
    \right| .
    \label{dsip}
\end{eqnarray}
Therefore
\begin{multline}
    \sqrt{1-\delta'}\int_{B(p,3\eps)}\left|d s_i^p(y)\right|dV(y) \leq\\
    \left(V(p,3\eps)
     \int_{B(p,3\eps)}|\nabla s(y)|^2dV(y)\right)^{\frac{1}{2}}
     + n \sum_{j=1}^n \|\nabla e_j^p\|_{2,3\eps}
     \|s_j^p\|_{2,3\eps} \nonumber
\end{multline}
so that we obtain by Lemma \ref{lmmDev} and by construction of
$e_j^p$
\begin{eqnarray}
    \sum_{i=1}^n \left( \int_{B(p,3\eps)}\left|d s_i^p(y)
    \right|dV(y)\right)^2 \leq
    c  \int_{B(p,3\eps)}|\nabla s(y)|^2dV(y)
     + c \widetilde{V}s(p) .\nonumber
\end{eqnarray}
We have then the following upper bound for (\ref{eqn8})
\begin{multline}
    \sum_{i=1}^n \left( \frac{1}{V'_{pq}}\int_{B'_{pq}}\left|
        \D s(p)_i - s_i^p(y)
    \right|dV(y)\right)^2  \label{eqn10} \\
    \leq c_K^2 c \left( \int_{B(p,3\eps)}|\nabla s(y)|^2dV(y)
     +  \left(\widetilde{V}s\right)(p)\right) .
\end{multline}
By the same kind of arguments as for (\ref{eqn8}) and using that
$\sum_{i,j=1}^n |A(q,p)_{ij}|^2$ is bounded above by a uniform
constant, we can bound (\ref{eqn7}) as follows
\begin{multline}
    \sum_{i=1}^n \left( \frac{1}{V'_{pq}}\int_{B'_{pq}}\left|
        \sum_{j=1}^n A(q,p)_{ij} \left(s_j^q(y)-\D s(q)_j \right)
    \right|dV(y)\right)^2 \label{eqn11} \leq \\
    c' \left( \int_{B(q,3\eps)}|\nabla s(y)|^2dV(y)
     + \left(\widetilde{V}s\right)(q)\right) .
\end{multline}
The last term (\ref{eqn9}) is then bounded by the following lemma
\begin{lmm}\label{lmm324}
    There exists a positive constant $c$ depending only on $m$, $n$, $k_1$,
    $k_2$, $\kappa$ and $\eps$ such that
    \begin{multline}
    \sum_{i=1}^n \left(\int_{B'_{pq}}\left|
    s_i^p(y)-\sum_{j=1}^n A(q,p)_{ij} s_j^q(y) \right|dV(y)\right)^2
    \leq \\
    c \left(\left(\widetilde{V}f\right)(p)+
    \left(\widetilde{V}f\right)(q)\right) .\nonumber
\end{multline}
\end{lmm}
\textbf{Proof}: see Appendix \ref{applmm2}. $\Box$

\medskip

Finally, (\ref{eqn10}), (\ref{eqn11}) and Lemma \ref{lmm324} imply
that
\begin{eqnarray}
    |D_A (\D s )(p,q)|^2 &\leq& c''\left(
    \int\limits_{B(p,3\eps)}|\nabla s(y)|^2dV(y) +
    \int\limits_{B(q,3\eps)}|\nabla s(y)|^2dV(y)\right)\nonumber\\
    & & + \
    c '' \left( \left(\widetilde{V}s\right)(p) +
    \left(\widetilde{V}s\right)(q)\right) .
    \nonumber
\end{eqnarray}
Taking the sum over $p$ and $q$ leads to
\begin{eqnarray}
    \|D_A(\D s)\|^2 & \leq &
    c''' \left(\|\nabla s\|^2 +
    \sum_{p\in X}\left(\widetilde{V}s\right)(p)\right) \label{eqn14}
\end{eqnarray}
so that (\ref{eqn14}),(\ref{eqn15}) imply (\ref{eqn13}). In order
to conclude the proof of point $ii)$ of this lemma, we have to
show that there exists $c >0$ such that
\begin{eqnarray}
    \sum_{p\in X}\left(\widetilde{V}s\right)(p) \leq c \|\nabla
    s\|^2 .\label{eqn16}
\end{eqnarray}
Fix $q \in X$, let $B = B(q,10 \eps)$, $V(B)$ its volume. Let
$(\cdot,\cdot)_B$ and $\|\cdot\|_B$ the $L^2$-inner product
respectively the $L^2$-norm on $E$ restricted to $B$. We are going
to show that there exists $c > 0$ such that
\begin{eqnarray}
    \left(\widetilde{V}s\right)(q)  \leq
    c \sum_{p \in B(q,3\eps)\cap X}\|\nabla s\|^2_{B(p,10\eps)} .\label{eqn18}
\end{eqnarray}
Then (\ref{eqn16}) is a direct consequence of (\ref{eqn18}). To
prove (\ref{eqn18}) we have to consider separately the cases $E$
is of complex (or quaternionic) rank one and $E$ is of harmonic
curvature.

\medskip

\textbf{Assume $\mathbf{E}$ is of rank one}. The proof of
(\ref{eqn18}) in this case is much easier than in the other case
as the potential involves only the first eigenvalue of the ball.
Recall that $\lambda_1(q) \leq \frac{\|\nabla s \|^2_B}{\| s
\|^2_B}$ for any non-zero $s$. Therefore and as $B(q,3\eps)
\subseteq B(p,10\eps)$ for any $p\in N(q)$
\begin{multline}
    \left(\widetilde{V}s\right)(q)
    \leq  \|s\|^2_{B(q,3\eps)}
    \sum_{p\in B(q,3\eps)\cap X}\frac{\|\nabla s \|^2_{B(p,10\eps)}}{\| s \|^2_{B(p,10\eps)}}
      \leq  \sum_{p \in B(q,3\eps)\cap X}\|\nabla s\|^2_{B(p,10\eps)}\nonumber
\end{multline}
and this concludes the first case.

\medskip

\textbf{Assume $\mathbf{E}$ is of harmonic curvature}. If $y \in
B$, write $s(y)$ as a sum of orthogonal sections (with respect to
$(\cdot,\cdot)_B$) $s(y) = \widetilde{s}(y) + s^{\bot}(y)$ with
$\widetilde{s}(y) = \sum_{j\leq \mu(q)} \frac{(s,e_j^q)_B}{ V(B)}
e_j^q(y)$. We have the following properties of the decomposition.
\setlength{\columnsep}{0pt} \setlength{\multicolsep}{0pt}
\begin{multicols}{2}
\begin{eqnarray}
    (s^{\bot},e_j^q)_B    & =&     0 \;, \; \; \forall j \leq \mu(q), \nonumber\\
    \|s\|^2_B       &=&
    \|s^{\bot}\|_B^2 + \|\widetilde{s}\|^2_B ,\nonumber\\
    \|\widetilde{s}\|_B^2 &=&
    \sum_{j\leq \mu(q)}\frac{(s,e_j^q)_B^2}{V(B)} ,\nonumber
\end{eqnarray}

\begin{eqnarray}
    (\nabla s^{\bot} , \nabla \widetilde{s})_B &=&   0 ,\nonumber\\
    \|\nabla s\|^2_B      & =&
    \|\nabla s^{\bot}\|_B^2 + \|\nabla \widetilde{s}\|^2_B ,
    \nonumber\\
        \|\nabla \widetilde{s}\|_B^2& =&
    \sum_{j\leq \mu(q)}\frac{(s,e_j^q)_B^2}{V(B)} \lambda_j(q) .\nonumber
\end{eqnarray}
\end{multicols}
Consider then two cases. First assume $\|s^{\bot}\|^2_B=0$. Then
$s(y) = \widetilde{s}(y)$ which means that if $y \in B(p,10\eps)$
\begin{eqnarray}
    s_j^q(y) &=&    \left\{ \begin{array}{ll}
                            0 & \text{if } j > \mu(q) ,\\
                            \frac{(s,e_j^q)_B}{V(B)} & \text{if }
                            j\leq \mu(q) .
                            \end{array}
                    \right.
\nonumber
\end{eqnarray}
Therefore
\begin{eqnarray}
    \left(\widetilde{V}s\right)(q) & = &
    \left( \sum_{j\leq \mu(q)}\lambda_j(q)\int_{B(q,3\eps)}|s_j^q|^2dV +
    \sum_{j > \mu(q)} \int_{B(q,3\eps)}|s_j^q|^2dV \right)\nonumber\\
    & = &
    V(q,3\eps)\sum_{j\leq \mu(q)}\frac{(s,e_j^q)_B^2}{V(B)^2}
    \lambda_j(q)\nonumber
    \ \leq \
    c \|\nabla \widetilde{s}\|_B^2 .\nonumber
\end{eqnarray}
Moreover as
$s^{\bot}$ is zero $\|\nabla \widetilde{s}\|_B^2 =\|\nabla
s\|_B^2$ and so in this case (\ref{eqn18}) is verified.

\medskip

For the second case, assume $\|s^{\bot}\|_B^2 \neq 0$. Then apply
max-min theorem to $s^{\bot}$ to obtain $ \lambda_{\mu(q)+1}(q)
\leq \frac{\|\nabla s^{\bot}\|^2_B}{\|s^{\bot}\|^2_B}$ and by
definition of $\mu(q)$ this implies that
\begin{eqnarray}
    \delta \|s^{\bot} \|_B^2 \leq\|\nabla s^{\bot}\|_B^2 .\label{eqn20}
\end{eqnarray}
Moreover, let us rewrite $s^{\bot}$ as follows, for $y \in
B(q,8\eps)$
\begin{eqnarray}
    s^{\bot}(y) & = &
    \sum_{j\leq \mu(q)} \left(s_j^q(y) -\frac{(s,e_j^q)_B}{V(B)} \right)e_j^q(y) +
    \sum_{j> \mu(q)}s_j^q(y)e_j^q(y) .\nonumber
\end{eqnarray}
As $\{e_j^q(y)\}$ is an almost orthonormal basis, we obtain for $y
\in B(q,8 \eps)$
\begin{eqnarray}
    \sum_{j\leq \mu(q)} \left|s_j^q(y) -\frac{(s,e_j^q)_B}{V(B)} \right|^2 +
    \sum_{j> \mu(q)}|s_j^q(y)|^2
    & \leq &
    (1-\delta')^{-1}|s^{\bot}(y)|^2 .
    \nonumber
\end{eqnarray}
In particular, this implies
\begin{eqnarray}
    \sum_{j > \mu(q)}\int_{B(q,3 \eps)} |s_j^q(y)|^2dV(y) &\leq& (1-\delta')^{-1}
    \|s^{\bot}\|_B^2 \label{eqn21}
\end{eqnarray}
and
\begin{multline}
    \sum_{j\leq \mu(q)}\lambda_j(q)
    \int_{B(q,3\eps)}|s_j^q(y)|^2dV(y)\leq \\
     2\sum_{j\leq \mu(q)}\lambda_j(q) \int_{B(q,3\eps)}
        \left|s_j^q(y)-\frac{(s,e_j^q)_B}{V(B)}\right|^2dV(y)
    +2\sum_{j\leq
    \mu(q)}\frac{(s,e_j^q)_B^2}{V(B)}\lambda_j(q) \\
    \leq \frac{2\delta}{1-\delta} \|s^{\bot}\|^2_B + 2 \|\nabla \widetilde{s}\|^2_B .
    \label{eqn22}
\end{multline}
Then (\ref{eqn21}) and (\ref{eqn22}) imply that
$\left(\widetilde{V} s\right)(q)  \leq c \left( \|s^{\bot}\|^2_B +
\|\nabla \widetilde{s}\|^2_B \right)$. Use (\ref{eqn20}) together
with this inequality to obtain (\ref{eqn18}) and therefore
(\ref{eqn16}). Finally (\ref{eqn13}) together with (\ref{eqn16})
imply $ii)$.

\medskip

To prove $iii)$ consider the following sum. By the work of Buser
(Lemma 5.1 in \cite{Buser}), there exists $c_B >0$ depending only
on $m$, $\kappa$ and $\eps$ such that
\begin{eqnarray}
    \sum_{i=1}^n \int\limits_{B(p,3\eps)}
    \left| \D s(p)_i - s_i^p(x) \right|^2dV(x) & \leq &
    c_B \sum_{i=1}^n \int\limits_{B(p,3\eps)}
    \left| d s_i^p(x)\right|^2dV(x) .
    \nonumber
\end{eqnarray}
Moreover, using (\ref{dsip}) we obtain
\begin{multline}
    \sum_{i=1}^n \int\limits_{B(p,3\eps)}
    \left| \D s(p)_i - s_i^p(x) \right|^2dV(x)  \leq \\
    \frac{2 n c_B}{1-\delta'} \left( \int_{B(p,3\eps)}|\nabla s(y)|^2dV(y)
     +  n \sum_{j=1}^n \|\nabla e_j^p\|^2_{\infty,3\eps}
     \|s_j^p(y)\|^2_{2,3\eps}\right) .\label{starstar}
\end{multline}
Therefore, from (\ref{starstar}) we obtain
\begin{multline}
    \left| \D s(p) \right|^2  \geq c\!\!\!\!
    \int\limits_{B(p,3 \eps)}\!\!\!\sum_{i=1}^n
    \left| (s_i^p(x)-\D s(p)_i) - s_i^p(x) \right|^2dV(x)\\
    \geq
   c \!\!\!\!\int\limits_{B(p,3 \eps)}\!\!\!\sum_{i=1}^n |s_i^p(x)|^2 dV(x)
    - 2 c\!\!\!\!\int\limits_{B(p,3 \eps)}\!\!\!\sum_{i=1}^n|s_i^p(x)| |\D s(p)_i - s_i^p(x)|dV(x)
    \\
    \geq c'\|s\|^2_{B(p,3\eps)}-
     c''\|s\|_{B(p,3\eps)}
     \left(\|\nabla s\|^2_{B(p,3\eps)} +
     \sum_{j=1}^n \|\nabla
     e_j^p\|^2_{\infty,3\eps}\|s_j^p\|^2_{2,3\eps}
     \right)^{\frac{1}{2}}\label{dsp2}
\end{multline}
\textbf{Assume $\mathbf{E}$ is of harmonic curvature} and combine
Lemma \ref{LocExten} and Lemma \ref{lmmDev} with (\ref{dsp2}) to
obtain
\begin{eqnarray}
    \left| \D s(p) \right|^2
    \geq c'\|s\|^2_{B(p,3\eps)} -
     c''\|s\|_{B(p,3\eps)}
     \left(\|\nabla s\|^2_{B(p,3\eps)} + \left(\widetilde{V}s\right)(p)
     \right)^{\frac{1}{2}} .\nonumber
\end{eqnarray}
Moreover, by (\ref{eqn18}) $\left(\widetilde{V}s\right)(p)$ is
bounded above by $c\!\!\!\! \sum\limits_{q \in B(p,3\eps)\cap
X}\|\nabla s\|^2_{B(q,10\eps)}$. Then, taking the sum over $p \in
X$ produces new $c'$, $c'' >0$ such that
\begin{eqnarray}
    \|\D s\|^2  & \geq &
                c'\|s\|^2-c''\|s\|\|\nabla s\|   .       \nonumber
\end{eqnarray}
Finally, if $\|\nabla s\|^2 \leq \Lambda'\|s\|^2$, we get $ \|\D
s\|^2  \geq \|s\|^2(c'-c''\sqrt{\Lambda'})$. Choose then
$\Lambda'$ suitably to conclude the proof of the proposition in
this case.

\medskip

\textbf{Assume $\mathbf{E}$ is of rank one}. If $\lambda_1(p)\leq
\delta$, by Lemma \ref{LocExten}, $\|\nabla
e_j^p\|^2_{\infty,3\eps} \leq c \lambda_1^s(p)$. If $\lambda_1(p)
> \delta$, by Lemma \ref{lmmDev} $\|\nabla
e_j^p\|^2_{\infty,3\eps} \leq c \leq c \delta^{-1} \lambda_1(p)$.
Therefore, (\ref{dsp2}) can be changed in (with new constants $c$,
$c'$, $c''$ )
\begin{eqnarray}
    \left| \D s(p) \right|^2 \geq
    \left\{
    \begin{array}{ll}
       \scriptstyle{ (c'-c\lambda_1^{\frac{s}{2}}(p))\|s\|^2_{B(p,3\eps)}-
        c''\|s\|_{B(p,3\eps)}\|\nabla s\|_{B(p,3\eps)}} &
        \text{if } \lambda_1(p) \leq \delta ,\\
        &  \\
       \scriptstyle{ c'\|s\|^2_{B(p,3\eps)}- c''\|s\|_{B(p,3\eps)}
     \|\nabla s\|_{B(p,10\eps)}}&
        \text{otherwise.}
    \end{array}
    \right.
    \nonumber
\end{eqnarray}
By choosing $\delta$ smaller, we can assume that if
$\lambda_{1}(p) \leq \delta$, $c'-c\lambda_1(p)^{\frac{s}{2}} \geq
c''' > 0$. This implies that (for any values of $\lambda_1(p)$)
\begin{eqnarray}
    \left| \D s(p) \right|^2 & \geq &
    c'''\|s\|^2_{B(p,3\eps)}-
    c''\|\nabla s\|_{B(p,3\eps)}\|s\|_{B(p,10\eps)} .\nonumber
\end{eqnarray}
Then, take the sum over $p \in X$ to obtain for $\|\nabla s \|
\leq \Lambda' \|s\|$
\begin{eqnarray}
    \| \D s\|^2 & \geq &
    c''' \|s\|^2 - c''\|\nabla s\| \|s\| \nonumber \ \geq \
    \|s\|^2(c'''-c''\sqrt{\Lambda'}) \nonumber
\end{eqnarray}
and conclude choosing $\Lambda'$ suitably. $\Box$

%%%%%%%%%%%%%%%%%%%%%%%%%%%%%%%%%%%%%%%%%%%%%%%%%%%%%%%%%%%%%%%%%%%%%%%%%%%%%%%%%%%%%%%%

\subsection{Upper bounds}

\begin{lmm}\label{lmmCVBUB}
    Let $m$, $n$, $k_1$, $k_2$, $\kappa$, $r_0$, $\eps$ be as before. Then
    there exist positive constants $c_3$ and $c_3'$ depending only
    on $m$, $n$, $k_1$, $k_2$, $\kappa$, $\eps$ so that for any vector bundle $E \in
    \E$ over any $M \in \M$, for any $X$ $\eps$-discretization of
    $E$ and for $\Delta_A +V$ constructed in Section \ref{Constr},
    we have
    \begin{enumerate}
        \item[i)] $\lambda_k(E) \leq c_3$, $\forall k \leq n |X|$,
        \item[ii)] $\lambda_k(X,A,V) \leq c'_3$, $\forall k \leq n |X|$.
    \end{enumerate}
\end{lmm}

\textbf{Proof}: $i)$ Let $p_i$ be a vertex of $X$ and consider
$f_i : M \to \R$ the first eigenfunction of the Dirichlet problem
for the ball centered at $p_i$ of radius $\frac{\eps}{2}$ extended
by zero. By Cheng's comparison theorem
$\frac{\|df_i\|^2}{\|f_i\|^2} \leq
\lambda_1\left(\frac{\eps}{2},\kappa\right)$ (where
$\lambda_1\left(\frac{\eps}{2},\kappa\right)$ denotes the first
non-zero eigenvalue of the Dirichlet problem on the ball of radius
$\frac{\eps}{2}$ in the simply connected space of constant
sectional curvature $-\kappa$ and of same dimension as $M$).
Define then the sections $\sigma_j^i (x) = f_i (x) e_j^{p_i}(x)$
for $1\leq i \leq |X|$, and $1 \leq j \leq n$. Then $\{
\sigma_{j}^{i} \ | \ 1 \leq i \leq |X| , 1 \leq j \leq n  \}$
spans a vector subset $W$ of $\Gamma(E)$ of dimension $n |X|$ as
$\{e_j^{p_i}\}_{j=1,\ldots,n}$ is an almost orthonormal frame.
Moreover
\begin{eqnarray}
    \nabla \sigma_j^i(x) & = & df_i(x)e_j^{p_i}(x) + f_i(x)\nabla
    e_j^{p_i}(x)\nonumber
\end{eqnarray}
hence by construction of $e_j^{p_i}$ and Lemma \ref{LocExten} and
Lemma \ref{lmmDev}, we have
\begin{eqnarray}
    \|\nabla \sigma_j^i\|^2 & \leq & c\left(\|df_i\|^2 +
    \|f_i\|^2\right)
    \nonumber
\end{eqnarray}
so that by definition of the $f_i$'s
\begin{eqnarray}
    \|\nabla \sigma_j^i\|^2 & \leq & c \|f_i\|^2
    \left( 1 + \lambda_1\left(\frac{\eps}{2},\kappa\right)
    \right) .
    \nonumber
\end{eqnarray}
By min-max theorem we get then
\begin{eqnarray}
    \lambda_k(E) & \leq
    c' \max\left\{ \frac{\sum_{i,j} a_{ij}^2\|\nabla \sigma_j^i\|^2}{\sum_{i,j} a_{ij}^2 \|\sigma_j^i\|^2}
    \right\}
    & \leq c' c \left( 1 + \lambda_1\left(\frac{\eps}{2},\kappa\right)\right) .\nonumber
\end{eqnarray}
This concludes the first part of the lemma.

\medskip

$ii)$ Let $f \in \FnX$. As $A(p,q)$ is a change of almost
orthonormal bases we have
\begin{eqnarray}
    \|D_A f\|^2 + (Vf,f) & = &
    \frac{1}{2}\sum_{p\in X}\sum_{q \in N(p)}|f(q) - A(p,q)f(p)|^2
    + \sum_{p\in X} (Vf)(p) \cdot f(p)\nonumber\\
    & \leq &
    c \sum_{p\in X}\sum_{q \in N(p)}\left( |f(p)|^2 + |f(q)|^2 \right)
    + \max\{\delta,1\}\|f\|^2\nonumber\\
    & \leq &
    \left(2 c \nu_X + \max\{\delta,1\}\right)\|f\|^2 .\nonumber
\end{eqnarray}
Therefore, $R(f) \leq 2 c \nu_X+ \max\{\delta,1\}$, $\forall f \in
\FnX \setminus\{0\}$ and this implies $\lambda_k(X,A,V) \leq 2 c
\nu_X + \max\{\delta,1\} $, $\forall k \leq n |X|$. $\Box$

%%%%%%%%%%%%%%%%%%%%%%%%%%%%%%%%%%%%%%%%%%%%%%%%%%%%%%%%%%%%%%%%%%%%%%%%%%%%%%%%%%%%%%%%

\subsection{Conclusion}\label{concl}

\textbf{Proof of Theorem \ref{main}}: by symmetry  of the results
concerning the smoothing and the discretizing, it suffices to
deduce $\lambda_k(E) \leq c \lambda_k(X,A,V)$. The proof proceeds
in two steps.

\medskip

First, assume that $k$ is such that $\lambda_k(X,A,V) \geq
\Lambda$, for $\Lambda$ given by Proposition \ref{lmmCVBS} iii).
Then, Lemma \ref{lmmCVBUB} $i)$ leads to $\lambda_k(E) \leq c_3
\Lambda^{-1} \lambda_k(X,A,V)$.This is the required inequality.

\medskip

Secondly, assume that $k$ is such that $\lambda_k(X,A,V) \leq
\Lambda$. Let $W_k$ be the $k$-dimensional vector subspace of
$\FnX$ spanned by $f_i : X \to \R^n$, $i=1,\ldots,k$,
$\lambda_i(X,A,V)$-eigenfunction of $\Delta_A$ chosen so that
$(f_i,f_j) = \delta_{ij} |X|$. By min-max theorem,
$\lambda_k(X,A,V) = \max\{R(f) : f \in W_k\setminus\{0\}\}$. Let
then $\S W_k$ be the vector subspace of $\Gamma(E)$ spanned by the
$\S f_i$'s i.e. $\S W_k =\langle\S f_1 , \ldots, \S f_k \rangle
=\{\S f \ | \ f \in W_k \setminus \{0\}\}$. As $\lambda_k(X,A,V)
\leq \Lambda$, for any non-zero function $f$ in $W_k$, we have
$\|D_A f\|^2 + \left(Vf,f\right) \leq \Lambda \|f\|^2$. Hence, by
Proposition \ref{lmmCVBS} $iii)$, for any $f$ in $W_k$, $\|\S
f\|^2 \geq c_2 \|f\|^2$ holds. In particular, $\S f$ is the zero
function if and only if $f$ is zero which means that $\S W_k$ is
$k$-dimensional. So we can apply min-max theorem to $\S W_k$ and
obtain
\begin{eqnarray}
    \lambda_k(E) \leq \max \{ R(\S f) \ | \ f \in W_k
    \setminus\{0\}\} .
    \nonumber
\end{eqnarray}
Moreover, by Proposition \ref{lmmCVBS} $ii)$ and $iii)$ we obtain
that $R(\S f) \leq \frac{c_1}{c_2} R(f)$ for any non-zero $f$ in
$W_k$, which leads to
\begin{eqnarray}
    \lambda_k(E) \leq \frac{c_1}{c_2} \max \{ R(f) \ | \ f \in W_k
    \setminus \{0\}\} = \frac{c_1}{c_2} \lambda_k(X,A,V) .\nonumber
\end{eqnarray}
This concludes the proof. $\Box$

\section{Estimation of the first non-zero eigenvalue for a flat vector bundle}\label{holonomy}

Let $(E^n, \nabla)$ be a flat Riemannian vector bundle with
irreducible holonomy over $M \in \M$. We recall the definition of
the constant related to the holonomy given by Ballmann, Br\"uning
and Carron in \cite{BBC}. If $c$ is a unit speed loop, denote by
$H_c$ its holonomy. Then there exists $\alpha > 0$ such that
$\forall x \in M$, $\forall v \in E_x$ there exists a smooth unit
speed loop $c_{x,v}$ of length less than two diameters of $M$ such
that
\begin{eqnarray}
    | H_{c_{x,v}} (v) - v | \geq \alpha |v| .\label{defholo}
\end{eqnarray}
The following theorem shows that if $E$ has significant holonomy,
then the first eigenvalue of $\overline{\Delta}$ can not be too
small. Conversely, if there exists $v$ in $E_x$ which has a small
holonomy, then the first eigenvalue is not too large.
\begin{thm}\label{thmHolo}
    Let $(E^n, \nabla)$ be a flat Riemannian vector bundle over $M
    \in \M$ with irreducible holonomy. Then there
    exist $c$, $c' > 0$ depending only on  $m$, $n$, $\kappa$,
    $r_0$ such that
    \begin{eqnarray}
        \lambda_1(E) \geq c' \frac{\alpha^2}{d(M)^2 c^{d(M)}}\nonumber
    \end{eqnarray}
    where $d(M)$ denotes the diameter of $M$.

    \medskip

    Moreover, if there exist $p_0 \in M$, $v_0 \in E_{p_0}$ and
    $\alpha'$ such that for any loop $c$ at $p_0$ of length less
    than $7 d(M)$, $|H_c(v_0)-v_0| \leq \alpha' |v_0|$ then,
    there exists $c'' >0$ depending only on $n$, $m$, $\kappa$ and
    $r_0$ such that
    \begin{eqnarray}
        \lambda_1(E) \leq c'' \alpha'^2 .\nonumber
    \end{eqnarray}
\end{thm}
The first part of the theorem is in fact due to Ballmann, Brüning
and Carron (see \cite{BBC}). We present here a more conceptual
proof that relies on the fact that the discrete magnetic Laplacian
associated to a discretization of a flat bundle is strongly
related to the holonomy of the vector bundle.

\medskip

\textbf{Proof}: let $\eps = \frac{1}{100}r_0$ and let $X$ be an
$\eps$-discretization of $E$. Then by Theorem \ref{main} there
exist $\Delta_A$ a discrete magnetic Laplacian and $c > 0$ such
that $\lambda_1(E) \geq c \lambda_1(X,A)$. So it suffices to prove
the statement for $\lambda_1(X,A)$. Let $f \in \FnX$ such that
$\Delta_A f = \lambda f$. Let $p_0 \in X$ and $v_0 = \sum_{i=1}^n
f_i(p_0) e_i^{p_0} \in E_{p_0}$. By (\ref{defholo}), there exists
a smooth unit speed loop $c_0 : [0,l] \to M$ at $p_0$ of length $l
\leq 2 d(M)$ and $|H_{c_0}(v_0)-v_0| \geq \alpha |v_0|$. Let $N
\in \N$ such that $ N \frac{\eps}{2} \leq l < (N+1)\frac{\eps}{2}$
and consider a partition of $[0,l]$, $0=t_0 < t_1 < \ldots <
t_{N-1} < t_N = l$ such that $\frac{\eps}{2} \leq t_j - t_{j-1}
\leq \eps $. By definition of $X$, $\forall j = 1 , \ldots ,N-1$,
$\exists p_j \in X$ such that $d(p_j, c_0(t_j)) < \eps$. Moreover,
let $p_N = p_0 \in X$. Note that $d(p_{j-1},p_j) < 3 \eps$.
Consider then the piecewise geodesic loop $\overline{c}_0$ at
$p_0$ passing through all $p_j$, $j = 1, \ldots N-1$ (i.e
$\overline{c}_0$ joins $p_{j-1}$ to $p_j$ via the minimizing
geodesic $p_{j-1} p_j$). Note that $\overline{c}_0$ is of length
less than $ 3 N \eps \leq 12 d(M)$. Moreover, as $E$ is flat, the
holonomy of $c_0$ is the same as the holonomy of $\overline{c}_0$.
More precisely, parallel translation from $c_0(t_{j-1})$ to
$c_0(t_j)$ along $c_0$ is the same as parallel translation along
minimizing geodesics from $c_0(t_{j-1})$ to $p_{j-1}$, then from
$p_{j-1}$ to $p_j$ and finally from $p_j$ to $c_0(t_j)$. Hence
$H_{c_0}(v) = H_{\overline{c}_0}(v)$ for any $v \in E_{p_0}$. So
that we obtain
\begin{eqnarray}
    |H_{\overline{c}_0}(v_0) - v_0| \geq \alpha |v_0| = \alpha |f(p_0)| .\nonumber
\end{eqnarray}
Consider then $v_j = \sum_{i=1}^n f_i(p_j) e_i^{p_j} \in E_{p_j}$.
By triangle inequality and as parallel transport is an isometry,
we obtain easily the following inequality
\begin{eqnarray}
    \alpha |f(p_0)| \leq \sum_{j = 1}^N |\tau_{p_j,p_{j-1}}v_{j-1}-v_{j} | .\nonumber
\end{eqnarray}
Moreover, by construction of $D_A$ we have
\begin{eqnarray}
    |\tau_{p_j,p_{j-1}}v_{j-1}-v_{j}| &=&
    \left |\sum_{i=1}^n f_i(p_{j-1})\tau_{p_j,p_{j-1}}e_i^{p_{j-1}} - \sum_{i=1}^n f_i(p_j)e_i^{p_j}\right|\nonumber\\
    &=&
    \left |\sum_{i=1}^n\left( \sum_{k=1}^n A(p_{j-1},p_j)_{ik} f_k(p_{j-1}) -  f_i(p_j)\right)e_i^{p_j}\right|\nonumber\\
    &=&
    |D_A f(p_{j-1},p_j)| .\nonumber
\end{eqnarray}
This implies that $ \alpha |f(p_0)| \leq |D_A f(p_0,p_1)| + \ldots
+ |D_A f(p_{N-1},p_N)|$. We have shown that for any $p_0 \in X$,
there exists a piecewise geodesic loop $\overline{c}_0 =
\{p_0,p_1,\ldots,p_N\}$ of length less than $12 d(M)$ such that
\begin{eqnarray}
    \alpha^2 |f(p_0)|^2 \leq 4\frac{d(M)}{\eps}
    \left(|D_A f(p_0,p_1)|^2 + \ldots + |D_A f(p_{N-1},p_N)|^2\right)
    \nonumber
\end{eqnarray}
and $d(p_{j-1},p_j) <3 \eps$. The goal is to apply this last
inequality to $\|f\|^2$. To that end, we need to find an upper
bound for the number of loops of the kind $\{p,q,\ldots,p\}$ that
can pass through $p \in X$ and $q \in N(p)$ and of length less
than $12 d(M)$. This upper bound on the length of the loop implies
that such a loop can pass through at most $P \leq 12
\frac{d(M)}{\eps}$ points of $X$. Therefore, there are at most
$\nu^{P-1}$ loops of the kind $\{p,q,\ldots,p\}$ and each of these
loops is suitable for $P$ points in $X$. Hence, we obtain
\begin{eqnarray}
    \alpha^2 \|f\|^2 &\leq&
    P\nu^{P-1}8\frac{d(M)}{\eps}\|D_A f\|^2\nonumber\\
    &\leq&
    72\frac{d(M)^2}{\eps^2}\nu^{12\frac{d(M)}{\eps}}\|D_A f\|^2 .\nonumber
\end{eqnarray}
This leads then to the conclusion of the first part $
\alpha^2\frac{\eps^2}{72 d(M)^2 \nu^{12\frac{d(M)}{\eps}}} \leq
\lambda$.

\medskip

To prove the second part of the theorem let $\eps=
\frac{1}{100}r_0$ and $X$ be an $\eps$-discretization of $E$ such
that $p_0 \in X$. Recall that $X$ is the set of vertices of a
finite connected graph $G$. Then construct a spanning tree $S$ of
$G$ (see \cite{Bo}, Section I.2) as follows. Let $X_i = \{p \in X
\ | \ d_G(p,p_0) = i\}$ where $d_G$ denotes the path metric on
$G$. Note that if $q$ is in $X_i$ then there exists $q'$ in
$X_{i-1}$ which is joined by an edge to $q$. Let then $S$ be the
subgraph of $G$ with vertices set $X$ and edges set $E(S) = \{q q'
\ | \ q \neq p_0\}$. We have constructed a spanning tree $S$ of
$G$.

\medskip

By construction of $S$, for any $p$ in $X$ there exists a unique
curve $c_p$ in $S$  joining $p$ to $p_0$ (i.e. $c_p$ is a
piecewise geodesic curve $\{p, \ldots ,p_0\}$ such that two
consecutive points of $X$ in $c_p$ are joined in $S$). Moreover
the length of such a $c_p$ is bounded above by $3 d(M)$. Now,
choose in $E_{p_0}$ an orthonormal basis
$\{e_1^{p_0},\ldots,e_n^{p_0}\}$ and define an orthonormal basis
$\B_p$ of $E_p$ by $\B_p = \{e_i^p =
\tau_{c_p}e_i^{p_0}\}_{i=1,\ldots,n}$, where $\tau_{c_p}$ denotes
parallel transport along $c_p$ from $p_0$ to $p$. Then $e_i^p(x) =
\tau_{x,p}e_i^p$ gives a local orthonormal frame made of parallel
sections. Hence, consider the discrete magnetic Laplacian
$\Delta_A$ associated to this choice of bases (constructed as in
Section \ref{Constr}) which satisfies $\lambda_1(E) \leq c
\lambda_1(X,A)$ by Theorem \ref{main}. So that it suffices to
prove the result for the first eigenvalue of $\Delta_A$. By
min-max theorem $\lambda_1(X,A) \leq R(f)$ for any non-zero
function on $X$. So consider $f : X \to \R^n$ defined by $f(p) =
\sum_{i=1}^n v_i e_i$ where the $v_i$'s are the coordinates of
$v_0$ in the basis $\B_{p_0}$. If $p$ and $q$ are neighboring
points in $X$ such that $d(p,p_0) \leq d(q,p_0)$ and $p \in c_q$,
then we have $\tau_{q,p}e_j^p = e_j^q$. Hence in this case
$A(p,q)_{ij} = \delta_{ij}$ and so $D_A f(p,q) = 0$. In the other
case i.e. if $p \in N(q)$, $d(p,p_0) \leq d(q,p_0)$ and $p$ is not
on $c_q$, consider the loop $c$ at $x_0$ going from $x_0$ to $p$
via $c_p$, from $p$ to $q$ via the minimizing geodesic $pq$ and
from $q$ to $x_0$ via $c^{-1}_q$. Then $c$ is of length less than
$7 d(M)$ and by assumption
\begin{eqnarray}
    |H_c(v_0) - v_0| \leq \alpha'|v_0| .\label{eqn4}
\end{eqnarray}
But, we have $H_c(v_0) = \tau^{-1}_{c_q}\tau_{q,p}\tau_{c_p}v_0$
and
\begin{eqnarray}
    \langle H_c(v_0),e_i^{p_0} \rangle = \left \langle
    \sum_{j=1}^n\tau_{q,p}e_j^p,e_i^q \right \rangle = \sum_{j=1}^n
    A(p,q)_{ij}v_j .
    \nonumber
\end{eqnarray}
Combining this last equality with (\ref{eqn4}) we obtain $\alpha'
|v_0| \geq |D_A f(p,q)|$. Finally, computing $\|D_A f\|^2$ leads
to
\begin{eqnarray}
    \|D_A f\|^2 \leq \frac{1}{2} \alpha'^2 \nu \|f\|^2 .
    \nonumber
\end{eqnarray}
So that the second part of the theorem follows. $\Box$

\appendix

\section{Appendix: technical tools}
The following lemma is a generalization of Lemma 11.1 in \cite{Li}
and a local version of Lemma 0.1 of \cite{PS2}.
    \begin{lmm}\label{lmmApp}
    Let $M \in \M$ and $u$ a non-negative function on the ball
    $B(p,R)$, with $R < \frac{1}{2}r_0$, such that $\Delta u \leq
    \alpha u + \beta$. Let $0 < \theta < 1$.
    Then there exist $c_1$, $c_2$, $c_3 >0$ (depending only on $m$,
    $n$, $\kappa$, $R$, $\alpha$ and $\beta$)
    and $0 < c(m) < s \leq 1$ such that
    \begin{eqnarray}
        \|u\|_{\infty,\theta R}\leq
        \left(\left(c_1 + c_2 \frac{1}{(1-\theta)^2}\right)^{c_3}  \|u\|_{2,R}
        \right)^s
        \nonumber
    \end{eqnarray}
    where $\|u\|_{\infty,\theta R} =\sup\{u(x) \ | \ x \in B(p, \theta R)\}$, and
    $\|u\|^q_{q,R} = \int_{B(p,R)} u^q(x)dV(x)$.
\end{lmm}
Note that, if $\beta = 0$ then $s=1$ (see \cite{Li}, Lemma 11.1).

\medskip

\textbf{Proof}: the proof combines the proof given in \cite{Li}
(Lemma 11.1) and Lemma 0.1 of \cite{PS2} . Let $u : B(p,R) \to
\R$, $u \geq 0$ such that $\Delta u \leq \alpha u + \beta$. Let $
\nu = \frac{m}{2}$ if $m\geq 3$ and $\nu = 2$ otherwise. Let $\mu$
be such that $\frac{1}{\mu} + \frac{1}{\nu} = 1$. For $0 < \rho <
\rho+\sigma < R$, let $\phi_{\rho, \sigma}$ be the Lipschitz
cut-off function depending only on the distance to $p$ given by
\begin{eqnarray}
    \phi_{\rho,\sigma}(r) =
    \phi(r) =\left\{\begin{array}{ll}
                    0 & \text{ on }
                    B(p,R)\setminus B(p,\rho+\sigma) ,\\
                    \frac{\rho+\sigma+r}{\sigma} & \text{ on }
                    B(p,\rho+\sigma)\setminus B(p,\rho) ,\\
                    1 & \text{ on }
                    B(p, \rho) .
                    \end{array}\right.
    \nonumber
\end{eqnarray}
Then for an arbitrary constant $a \geq 1$, we have
\begin{eqnarray}
    \|u^{2a}\|_{\mu,\rho} \leq \|\phi u^a\|^2_{2 \mu} .
    \nonumber
\end{eqnarray}
As the injectivity radius of $M$ is bounded below ($Inj(M)\geq r_0
>0)$ and the Ricci curvature too ($Ricci(M,g)\geq -(m-1)\kappa g$)
Sobolev embeddings for complete manifolds are valid  and we can
apply the Sobolev inequalities to $\|\phi u^a\|^2_{2\mu}$ (see
\cite{He1}, lemma 3.3). More precisely, there exists a constant
$c_s
> 0$ depending only on $m$, $\kappa$ and $r_0$ such that
\begin{eqnarray}
    \|\phi u^a\|^2_{2\mu} \leq c_s \|d(\phi u^a)\|_2^2 .
    \nonumber
\end{eqnarray}
Replacing $c_s$ by $C R^2$, we can rewrite the inequality as
\begin{eqnarray}
    \|\phi u^a\|^2_{2\mu} \leq C R^2 \|d(\phi u^a)\|_2^2 .
    \nonumber
\end{eqnarray}
Therefore,
\begin{eqnarray}
    \|u^{2a}\|_{\mu,\rho} \leq C R^2 \|d(\phi u^a)\|_2^2 .
    \nonumber
\end{eqnarray}
However
\begin{eqnarray}
    \int_{M} |d(\phi u^a)|^2dV \leq
    a \int_{M}\phi^2 u^{2a-1}\Delta u dV + \int_{M}|d\phi|^2 u^{2a}dV
    \nonumber
\end{eqnarray}
(see \cite{Li}, p.81). Hence using the assumption on $\Delta u$
and $u \geq 0$ we obtain
\begin{eqnarray}
    \|u^{2a}\|_{\mu,\rho} \leq    C R^2\left( a \alpha\int_M\phi^2u^{2a}dV
                                +  a \beta  \int_M \phi^2 u^{2a-1}dV
                                +  \int_M |d\phi|^2u^{2a}dV\right)
    \nonumber
\end{eqnarray}
and by construction of $\phi$, we obtain
\begin{multline}
    \|u^{2a}\|_{\mu,\rho}  \leq
    CR^2\left( a \alpha + \frac{1}{\sigma^2}\right)\int_{B(p,\rho+\sigma)}u^{2a}dV
    \nonumber
    +
    CR^2 a \beta \int_{B(p,\rho+\sigma)}u^{2a-1}dV \nonumber\\
     \leq
    CR^2\left(a \alpha + \frac{1}{\sigma^2}\right)\|u\|^{2a}_{2a,\rho+\sigma}
    \nonumber
    +
    CR^2 a \beta
    V(p,\rho+\sigma)^{\frac{1}{2a}}\|u\|^{2a-1}_{2a,\rho+\sigma} .
    \nonumber
\end{multline}
Finally, we have shown that for any $a \geq 1$, $0 < \rho < \rho+
\sigma < R$, we have
\begin{eqnarray}
    \|u\|^{2a}_{2a\mu,\rho} & \leq &
    CR^2\left(a\alpha+\frac{1}{\sigma^2}\right)\|u\|^{2a}_{2a,\rho+\sigma}+
    CR^2 a \beta
    V(p,\rho+\sigma)^{\frac{1}{2a}}\|u\|^{2a-1}_{2a,\rho+\sigma} .
    \nonumber
\end{eqnarray}
This was the first step of the proof. Now, we will proceed with a
Moser iteration. To that aim, let
\begin{eqnarray}
    a_0 = 1 ,\ a_1 = \frac{m}{m-2} =\mu ,& \ldots ,& a_i=\mu^{i}, \ldots \nonumber\\
    \sigma_0 =\frac{1-\theta}{2}R,\ \sigma_1 =\frac{1-\theta}{4}R, & \ldots,&
    \sigma_i = \frac{1-\theta}{2^{i+1}}R,\ldots\nonumber\\
    \rho_0=R-\sigma_0 ,\ \rho_1 = R-\sigma_0-\sigma_1 ,&\ldots,&
    \rho_i = R-\sum_{j=0}^i\sigma_j, \ldots \nonumber
\end{eqnarray}
and $\rho_{-1}=R$. Observe that $\rho_i>\theta R$ for any $i$ and
$\rho_i \to \theta R$ as $i \to \infty$. Moreover, for any $A_i$,
$B_i >0$
\begin{multline}
    (A_i+B_i)\min\{\|u\|_{2a_i,\rho_i+\sigma_i}^{2a_i},
                    \|u\|_{2a_i,\rho_i+\sigma_i}^{2a_{i-1}}\}
    \leq \\ A_i \|u\|_{2a_i, \rho_i+\sigma_i}^{2a_i}+
                    B_i\|u\|_{2a_i, \rho_i+\sigma_i}^{2a_{i-1}}
     \leq  (A_i+B_i)\|u\|_{2a_i,\rho_i+\sigma_i}^{b_i}\nonumber
\end{multline}
where $b_i$ is suitably chosen ($b_i \in \{2a_{i-1},2a_i\}$). Now
replace above $a$, $\rho$, $\sigma$ by $a_i$ respectively
$\rho_i$, $\sigma_i$ to obtain
\begin{eqnarray}
    \|u\|_{2a_{i+1},\rho_i} \leq
    \left(CR^2\left(a_i \alpha + \frac{1}{\sigma_i^2} +
    a_i \beta V(p,\rho_{i-1})^{\frac{1}{2a_i}}
    \right)\right)^{\frac{1}{2a_i}}\|u\|^{\frac{b_i}{2a_i}}_{2a_i,\rho_{i-1}}
    .
    \nonumber
\end{eqnarray}
Then iterate this inequality to obtain (using Bishop-Gromov
comparison theorem, Croke's inequality and $a_i \geq 1$)
\begin{eqnarray}
    \|u\|_{\infty,\theta R} \leq c\left(
    \prod_{i=0}^{\infty}\left(
        CR^2a_i\left(\alpha + c'\beta\right)+C\frac{R^2}{\sigma_i^2}
    \right)^{\frac{1}{2a_i}} \|u\|^{\frac{b_0}{2}}_{2,R}
    \right)^{\prod_{j=1}^{\infty}\frac{b_j}{2a_j}} .
    \nonumber
\end{eqnarray}
By the same argument as in \cite{PS2},
$\prod_{j=0}^{\infty}\frac{b_j}{2a_j}$ converges to $s \in
[e^{-(n-2)\frac{\ln(2)}{2}},1]$. It remains then to show that
$\prod_{i=0}^{\infty}\left(
        CR^2a_i\left(\alpha + c'\beta\right)+C\frac{R^2}{\sigma_i^2}
    \right)^{\frac{1}{2a_i}}$ converges too.

But we have that $\prod_{i=0}^{\infty}B^{\mu^{-i}} =
B^{\frac{\mu}{\mu-1}}$ (as $\mu >1$) and
$\sum_{i=0}^{\infty}i\mu^{-i}$ is finite, therefore
\begin{multline}
    \prod_{i=0}^{\infty}\left(
        CR^2\mu^i\left(\alpha + c'\beta\right)+4C\frac{4^i}{(1-\theta)^2}
    \right)^{\frac{1}{2\mu^i}}\leq \\
    \prod_{i=0}^{\infty}\max\{\mu,4\}^{\frac{i}{2\mu^i}}\left(
        CR^2\left(\alpha + c'\beta\right)+C\frac{4}{(1-\theta)^2}
    \right)^{\frac{1}{2\mu^i}}\\
    \leq c(\mu)\left(CR^2\left(\alpha +
    c'\beta\right)+C\frac{4}{(1-\theta)^2}
    \right)^{\frac{1}{2}\frac{\mu}{\mu-1}} .\nonumber
\end{multline}
This implies the claim. $\Box$

%%%%%%%%%%%%%%%%%%%%%%%%%%%%%%%%%%%%%%%%%%%%%%%%%%%%%%%%%%%%%%%%%%%%%%%%%%%%%%%%%%%%%%%
\subsection{Proof of Lemma \ref{lmm320}} \label{applmm1}

The proof differs according to the assumptions made on $E$.

\medskip

\textbf{Assume $\mathbf{E}$ is of harmonic curvature}. By Remark
\ref{rmk321} and Remark \ref{rmk432}, we have
\begin{multline}
    \sum_{i=1}^n f_i(p)e_i^p(x) - \sum_{i=1}^n f_i(q)e_i^q(x) =
    \nonumber\\
    \sum\limits_{i=1}^n f_i(p)\left(e_i^p(x)-\tau_{x,p}e_i^p(p)\right)
    + D_A f(q,p)_i \tau_{x,p}e_i^p(p) + f_i(q)\left(\tau_{x,p}e_i^q(p)-e_i^q(x)\right).\nonumber
\end{multline}
By Lemma \ref{LocExten} and as $d^* R^E =0$,
$|e_i^p(x)-\tau_{x,p}e_i^p|^2 \leq c \lambda_i(p)$ for $1 \leq i
\leq \mu(p)$ and $|\tau_{x,p}e_i^q(p)-e_i^q(x)|^2 \leq c
\lambda_i(q)$ for $1 \leq i \leq \mu(q)$. Moreover if $\mu(q)<
i\leq n$, $|\tau_{x,p}e_i^q(p)-e_i^q(x)|^2 \leq 4$. Therefore
\begin{multline}
    \left|\sum_{i=1}^n f_i(p)e_i^p(x) - \sum_{i=1}^n f_i(q)e_i^q(x)\right|^2
    \leq \\
    c' \left(\left|D_A f(q,p)\right|^2 + \left(Vf\right)(p)\cdot f(p) + \left(Vf\right)(q)\cdot f(q)\right)
    \nonumber
\end{multline}
which implies the lemma in this case.

\medskip

\textbf{Assume $\mathbf{E}$ is of rank one}, then
\begin{multline}
    \sum_{i=1}^n f_i(p)e_i^p(x) - \sum_{i=1}^n f_i(q)e_i^q(x) = \\
    \sum_{i=1}^n D_A f(q,p)_i e_i^p(x) + \sum_{j=1}^n
    f_j(q)\sum_{i=1}^n e_i^p(x)\left(A(q,p)_{ij}-a(q,p)_{ij}(x)
    \right).
    \nonumber
\end{multline}
By definition of $A(q,p)_{ij}$ and by the work of Buser (Lemma 5.1
in \cite{Buser}) there exists $c_B >0$ depending only on $m$,
$\kappa$ and $\eps$ such that
\begin{eqnarray}
\int_{B_{pq}}\left|A(q,p)_{ij}-a(q,p)_{ij}(x) \right|^2 dV(x)\leq
c_B \int_{B_{pq}}\left|d a(q,p)_{ij}(x) \right|^2dV(x) .\nonumber
\end{eqnarray}
Moreover
\begin{multline}
    (1-\delta')\sum\limits_{i=1}^n\left|d a(q,p)_{ij}(x) \right|^2\leq\\
     \left|\sum\limits_{i=1}^n d a(q,p)_{ij}(x) e_i^p(x) \right|^2
    =  \left|\nabla e_j^q(x)-\sum\limits_{i=1}^n a(q,p)_{ij}(x)
    \nabla e_i^p(x) \right|^2 \\
    \leq  c \left(\left|\nabla e_j^q(x)\right|^2 + \sum\limits_{i=1}^n
    \left| \nabla e_i^p(x) \right|^2\right).\nonumber
\end{multline}
As the bundle is of rank one, $\lambda_1(p) = \ldots =
\lambda_n(p)$. Therefore $\lambda_1(p)\leq \delta$ implies
$\int\limits_{B(p,10\eps)}\left| \nabla e_i^p(x) \right|^2dV(x)
\leq c \lambda_1(p)$. Otherwise $\int\limits_{B(p,10\eps)}\left|
\nabla e_i^p(x) \right|^2dV(x) \leq c \leq c
\delta^{-1}\lambda_1(p)$ by Lemma \ref{lmmDev}, which implies
\begin{eqnarray}
\int_{B_{pq}}\left|A(q,p)_{ij}-a(q,p)_{ij}(x) \right|^2 dV(x)\leq
c' \left(\lambda_1(p)+\lambda_1(q)\right) .\label{aqpij}
\end{eqnarray}
Hence
\begin{multline}
    \int_{B(q,\eps)}\left|\sum_{i=1}^n f_i(p)e_i^p(x) - \sum_{i=1}^n f_i(q)e_i^q(x)\right|^2 \leq \\
    c'' \left( \left|D_A f(q,p)\right|^2 + \left| f(q)\right|^2\left(
    \lambda_1(p) + \lambda_1(q)\right)\right) .\nonumber
\end{multline}
This concludes the proof of Lemma \ref{lmm320}. $\Box$
%%%%%%%%%%%%%%%%%%%%%%%%%%%%%%%%%%%%%%%%%%%%%%%%%%%%%%%%%%%%%%%%%%%%%%%%%%%%%%%%%%%%%%
\subsection{Proof of Lemma \ref{lmm324}}\label{applmm2}

The proof differs according to the assumptions made on $E$.

\medskip

\textbf{Assume $\mathbf{E}$ is of harmonic curvature}. As
$\{\tau_{y,p}e_i^p(p)\}_{i=1}^n$ is an almost orthonormal basis
and by Remark \ref{rmk432}
\begin{multline}
    \sum_{i=1}^n
    \left| s_i^p(y)-\sum_{j=1}^n A(q,p)_{ij} s_j^q(y) \right|^2 \leq\\
    (1-\delta')^{-1} \left| \sum_{i=1}^n s_i^p(y)\left(\tau_{y,p}e_i^p(p) -e_i^p(y)\right)
    + \sum_{i=1}^n s_i^q(y)\left(e_i^q(y)
    -\tau_{y,p}e_i^q(p)\right)\right|^2 .
    \nonumber
\end{multline}
Integrate then over $B'_{pq}$ and apply Lemma \ref{LocExten} to
obtain
\begin{eqnarray}
    \sum_{i=1}^n \int_{B'_{pq}}
    \left|s_i^p(y)-\sum_{j=1}^n A(q,p)_{ij} s_j^q(y) \right|^2dV(y)
    \leq c \left( \left(\widetilde{V}s\right)(p) +
    \left(\widetilde{V}s\right)(q)\right).
    \nonumber
\end{eqnarray}

\medskip

\textbf{Assume $\mathbf{E}$ is of rank one}. Recall that
$s_i^p(y)= \sum_{j=1}^n a(q,p)_{ij}(y)s_j^q(y)$. Hence
$$s_i^p(y)-\sum_{j=1}^n A(q,p)_{ij} s_j^q(y) =
\sum_{j=1}^n\left(a(q,p)_{ij}(y)- A(q,p)_{ij}\right)s_j^q(y).$$
Therefore
\begin{multline}
    \int_{B'_{pq}}
    \left|s_i^p(y)-\sum_{j=1}^n A(q,p)_{ij} s_j^q(y) \right|dV(y)
    \leq \\
    \|s\|_{2,3\eps}\sum_{j=1}^n\left( \int_{B'_{pq}}
    \left|a(q,p)_{ij}(y)- A(q,p)_{ij}
    \right|^2\right)^{\frac{1}{2}}dV(y).
    \nonumber
\end{multline}
Finally, as $B'_{pq} \subset B_{pq}$, inequality (\ref{aqpij})
implies
\begin{eqnarray}
    \sum_{i=1}^n\left( \int_{B'_{pq}}
    \left|s_i^p(y)-\sum_{j=1}^n A(q,p)_{ij} s_j^q(y)
    \right|dV(y)\right)^2
    \leq c \left(\widetilde{V}s\right)(p)\nonumber
\end{eqnarray}
and this concludes the proof of Lemma \ref{lmm324}. $\Box$

{\small Tatiana Mantuano

Universit\'e de Neuch\^atel

Institut de Math\'ematiques

rue Emile-Argand 11

2009 Neuch\^atel

Switzerland

e-mail : Tatiana.Mantuano@unine.ch }

\end{document}